\documentclass[a4paper,preprint,12pt]{elsarticle}
\usepackage[english]{babel}
\usepackage[latin1]{inputenc}
\usepackage{amssymb}
\usepackage{amsmath}         
\usepackage{amsthm}             
\usepackage{pb-diagram}              
\usepackage{epstopdf}      
\usepackage{graphicx}             
 \usepackage{epsfig}     
\usepackage{color}         
\usepackage{fancyhdr} 
 \usepackage{natbib}      
\usepackage{rotating}    
\usepackage[T1]{fontenc}        
\usepackage{url}                
\usepackage{mathtools}                 

\usepackage{vmargin}            
\usepackage[title]{appendix} 

\usepackage{lineno}
\usepackage{setspace}                      
\biboptions{square} 
\usepackage[colorlinks=true,linkcolor= blue,urlcolor=blue,citecolor=blue]{hyperref} 
\usepackage[ruled,vlined]{algorithm2e} 

\newcommand{\hh}{h}

\newcommand{\bu}{\bullet}
 \newcommand{\bo}[1]{\mathbf{#1}}  
\def\indic{\hbox{1\kern-.24em\hbox{I}}}      

\newtheorem{assump}{Assumption}{\bf}{\it}

\newcommand{\var}{\mathbb{V}}  
\newcommand{\esp}{\mathbb{E}}

\newcommand{\M}{f}

\newcommand{\T}{T}    
\newcommand{\N}{\mathbb{N}}

\newcommand{\norme}[1]{\left|\left| #1 \right|\right|_{2}}    
\newcommand{\norml}[1]{\left|\left| #1 \right|\right|_{1}}

\newcommand{\X}{X}

\newcommand{\R}{\mathbb{R}}

{\bf}{\it} 
{\bf}{\it}  

\newtheorem{theorem}{Theorem}{\bf}{\it}     
\newtheorem{lemma}{Lemma}{\bf}{\it}          
\newtheorem{rem}{Remark}{\bf}{\it} 
\newtheorem{corollary}{Corollary}{\bf}{\it} 

\catcode`\"=\active
\catcode`\"=\active
\def"{\og\ignorespaces}
\def"{{\fg}}

\makeatletter
\def\ps@pprintTitle{%
  \let\@oddhead\@empty
  \let\@evenhead\@empty
  \def\@oddfoot{\reset@font\hfil\thepage\hfil}
  \let\@evenfoot\@oddfoot
}
\makeatother

\begin{document}
\begin{frontmatter}           
\title{Dimension-free estimators of gradients of functions with(out) non-independent variables}
\author[a,b]{Matieyendou Lamboni\footnote{Corresponding author: matieyendou.lamboni[at]gmail.com or [at]univ-guyane.fr; Dec 24, 2025}}               
\address[a]{University of Guyane, Department DFR-ST, 97346 Cayenne, French Guiana, France}
\address[b]{228-UMR Espace-Dev, University of Guyane, University of R\'eunion, IRD, University of Montpellier, France}                                                    
                                                                                          
\begin{abstract}   
This study proposes a unified stochastic framework for approximating and computing the gradient of every smooth function evaluated at non-independent variables, using $\ell_p$-spherical distributions on $\R^d$ with $d, p\geq 1$. The upper-bounds of the bias of the gradient surrogates  do not suffer from the curse of dimensionality for any $p\geq 1$. Also, the mean squared errors (MSEs) of the gradient estimators are bounded by  $K_0 N^{-1} d$ for any $p \in [1, 2]$, and by  $K_1 N^{-1} d^{2/p}$ when $2 \leq p \ll d$ with $N$ the sample size and $K_0, K_1$ some constants. Taking $\max\left\{2, \log(d) \right\} < p \ll d$ allows for  achieving dimension-free upper-bounds of MSEs. In the case where $d\ll p< +\infty$, the upper-bound  $K_2 N^{-1} d^{2-2/p}/ (d+2)^2$ is reached with $K_2$ a constant. Such results lead to dimension-free MSEs of the proposed estimators, which boil down to estimators of the traditional gradient when the variables are independent. Numerical comparisons show the efficiency of the proposed approach.      
\begin{keyword}         
Dependent variables \sep Gradients \sep High-dimensional models \sep  \sep Optimal estimators \sep Tensor metric of non-independent variables\\      
\textbf{AMS}: 60H25,  49Qxx. 90C25, 90C30, 90C56, 68Q25, 68W20, 65Y20       
\end{keyword}           	       
\end{abstract}           
\end{frontmatter}                               
 \setpagewiselinenumbers 
 \modulolinenumbers[1]        
     
 
      
\section{Introduction}        
Non-independent variables are often described by their covariance matrices, distribution functions, copulas, weighted distributions (see e.g., \cite{rosenblatt52,nataf62,joe14,mcneil15,navarro06,sklar59,durante22,lamboni23mcap}), and dependency models provide explicit functions that link these variables together by means of additional independent variables (\cite{skorohod76,lamboni21,lamboni22b,lamboni23mcap,lamboni24uq}). Functions evaluated at non-independent input variables are widely encountered in different scientific fields, and analyzing such functions requires being able to determine and to compute the non-Euclidean (or dependent) gradients, that is, the gradients that account for the  dependencies among the inputs.\\  
 
Dependent gradients of functions evaluate at non-independent variables (i.e., $grad\M$) have been proposed in \cite{lamboni23axioms,lamboni24axioms}, and it is shown that such gradients boil down to the traditional (or Euclidean) gradient (i.e., $\nabla \M$) when the inputs are independent. Indeed, $grad\M =G^{-1} \nabla \M$ with $G$ the tensor metric of non-independent variables, which is known or can be computed using the distribution elements of such variables. Computing $grad\M$ or equivalently $\nabla\M$ in  high-dimensional settings and for time-demanding models using a number of model evaluations that is less than the dimensionality $d$ requires relying on the Monte-Carlo approach (\cite{spall92,spall00,ancell07,pradlwarter07,patelli10}). The Monte-Carlo (MC) approach is a consequence of the Stokes theorem, which claims that the expectation of a function evaluated at a randomized point around $\bo{x} \in \R^d$ is the gradient of a certain function. Instances of its application are randomized approximations of $\nabla\M$, derivative-free methods, zeroth-order stochastic optimizations, online learning algorithms and optimizations, (see \cite{nemirovsky83,agarwal10, bach16,shamir17,akhavan20,akhavan22,akhavan24,lamboni24axioms,berahas22,ma25} and references therein). Moreover, the MC approach is relevant for applications in which the computations of the gradients are impossible (\cite{akhavan20}). \\      
    
Different stochastic surrogates of gradients are available in the literature. Among others, there are approximations that rely on i) randomized kernels and random vectors that are uniformly distributed over the unit ball or on the unit sphere (\cite{polyak90,agarwal10,bach16,shamir17,akhavan20,akhavan22,akhavan24}), and ii) random vectors only (\cite{agarwal10,shamir17,berahas22,lamboni24axioms}). An overview of different methods for the approximations and estimations of the traditional gradient is available in \cite{berahas22,gasnikov23}, and an unbiased formulation of the gradient is provided in \cite{ma25}. Mean squared errors (MSEs) are often used to assess the qualities of the gradient estimators based on the sample size $N \geq 1$. While such estimators are unbiased and consistent estimators of the corresponding gradient surrogates, they are biased estimators of gradients in general. In zeroth-order stochastic optimizations or online learning algorithms and optimizations, the sample size is implicitly fixed at $N=1$. Thus, it is worth noting that the MSEs of gradient estimators in the context of i) direct computations of gradients,  and ii) optimization problems are equal up to a factor, depending on $N^{-1}$. For finite differences methods (FDMs), one may also think that $N=1$ always.\\    
                         
 For exact-evaluated functions (i.e., noiseless cases) and smooth functions, the upper-bounds of the biases provided in \cite{polyak90,bach16,akhavan20} depend on the dimensionality in general, while dimension-free upper-bounds of the biases are obtained in \cite{berahas22,akhavan24,gasnikov23} for $\nabla\M$ and in \cite{lamboni23axioms} for $grad\M$.          
Regarding MSEs, the convergence rate of the form $\left(d^2 N^{-1}\right)^r$ is reached in \cite{bach16} with $r>0$, while $d^2N^{-1}$ is obtained in \cite{lamboni23axioms} for any differentiable function (see also, \cite{duchi15} for the upper-bound of the variance of their estimator). Using random vectors that are uniformly distributed on the unit $\ell_2$-sphere (i.e., sphere-based MC), the upper-bound of the variance of the gradient estimator proposed in \cite{shamir17} is bounded by $d N^{-1}$ up to a constant. However, the bias analysis is not available in in that paper, making difficult to derive the corresponding rate of convergence. The works in \cite{berahas22,scheinberg22} complement such a study, leading to the rate $d N^{-1}$ for the unit $\ell_2$-sphere-based MC approach as well as for FDMs. Also, dimension-free upper-bounds of biases and the MSE-based rate of the form $d N^{-1}$ are obtained in \cite{akhavan24,gasnikov23} using either the unit $\ell_2$-sphere or $\ell_1$-sphere-based MC approach in addition to randomized kernels. Recall that the aforementioned results correspond mainly to smooth and noiseless functions with the possibility of being evaluated twice in the case of zeroth-order stochastic optimizations. Such assumptions remain the basic ones for this paper. \\                         
	      
So far, the MSE-based rates of convergence of the aforementioned methods and others depend on the dimensionality. Such rates seem to be in contradiction with the consistent computation of $\nabla\M$ using $N \ll d$, as suggested in \cite{patelli10} and used in zeroth-order stochastic optimizations ($N=1$). Likewise, queries about the  theoretical  and numerical advantage of randomized schemes over the traditional FDMs are discussed in \cite{scheinberg22}. Significant advantages of randomized schemes over FDMs have been expected since the seminar works in \cite{spall92,spall00}, and this paper addresses such an open problem using the $\ell_p$-spherical distributions or random vectors that are uniformly distributed over the $\ell_p$-ball with $p\geq 1$.\\        
    
In this paper, a unified stochastic framework for approximating the gradients of functions evaluated at non-independent variables has been proposed, including the traditional gradients. It relies on a set of $L \geq 1$ constraints (leading to $L$-point-based surrogates) and a class of $\ell_p$-spherical distributions or uniform distributions over the $p$-balls. The $\ell_p$-sphere-based surrogates of the dependent gradients generalize those proposed in \cite{duchi15,shamir17,berahas22,lamboni24axioms}, as the dependent gradients boil down to the traditional gradients in the case of independent variables. The convergence analysis of the proposed gradient estimators shows that the biases and the MSEs do not suffer from the curse of  dimensionality by properly choosing the values of $p$. Table \ref{tab:comp} reports a summary of such an analysis compared to prior, best-known results.  From this table, the proposed approach outperforms the best known results and allows for breaking down the course of dimensionality. Also, the derived results are still valid for zeroth-order stochastic optimizations with one ($L=1$) or two ($L=2$)-point-based exact evaluations of functions.  
 
\begin{table}[htbp]  
\begin{center}   
\begin{tabular}{lccccc}       
\hline      
\hline          
	Methods		& Biases &   Rates & Costs & Types of gradients \\       			
\hline	
 FDM  (\cite{berahas22})  &  $\sqrt{d}Mh$   & d  &  d+1, 2d  & $\nabla\M$ \\
 MC (\cite{shamir17}) & $N\!A$  &   $d N^{-1}+ N\!A^2$  &   $2d$  &   $\nabla\M$  \\ 
 MC (\cite{berahas22}) &  $M h$   &  $dN^{-1}$  &   d+1, 2d  &  $\nabla\M$    \\ 
 MC (\cite{akhavan24}) &  $4\sqrt{2} M h$   &  $dN^{-1}$    &   2d  &    $\nabla\M$    \\ 
MC (\cite{lamboni24axioms}) &  $M h$  &   $d^2 N^{-1}$  &  $Ld^2$  &   $grad\M,\, \nabla\M$   \\ 
This paper, $p \in [1, 2]$ &  $M h$  &   $d N^{-1}$  &  $Ld$  &   $grad\M,\, \nabla\M$   \\ 
This paper, $2\leq p \ll d$ &  $M h$   &  $d^{2/p} N^{-1}$  &   $Ld^{2/p}$  &  $grad\M,\, \nabla\M$ \\
This paper, $d \ll p$ &  $M h$  &   $\frac{d^{2-2/p}N^{-1}}{(d+2)^2}$  &   $L d^{-2/p}$  &   $grad\M,\, \nabla\M$ \\   
\hline                                                         
\hline                                   
    \end{tabular}                           
  \end{center}    
  \caption{Summary of upper-bounds of biases and rates of convergence achieved by algorithms for the gradient computations under the assumption of differentiable functions without noises. The quantity $M$ stands for either the first-order ($M_1$-Lipschitz) or second-order ($M_2$) H\"older constant, and $h$ is the bandwidth. Costs stand for the crude, minimal and total number of evaluations required.} 
			 	\label{tab:comp}           
\end{table}                                      
                   
Formulations of dependent gradients and surrogates of such gradients are provided in Section \ref{sec:prel} by making usd of i) a set of constraints, and ii) the Monte-Carlo approach based on a wide class of random vectors, such as  independent random variables	that are symmetrically distributed about zero; $\ell_p$-spherically distributed random vectors. Each surrogate of the dependent gradient is followed by its order of approximation.  Section \ref{sec:dfexpderi} examines the convergence analysis of the proposed surrogates of gradients by distinguishing the bias analysis, the MSEs and the rates of convergence. 
 A number of numerical comparisons is considered so as to assess the numerical efficiency of our approach (see Section \ref{sec:fdm}). Section~\ref{sec:con} concludes this work.   
                                                  
\section{Generalized surrogates of dependent gradients} \label{sec:prel}    
This section deals with basic definitions of gradients of functions with or without dependent variables;  new random vectors that will be needed in the formulation of randomized surrogates of gradients;  new insight into the gradient surrogates based on $L$-point evaluations of functions and their corresponding estimators. 

\subsection{Notation}  
Denote with  $||\cdot||_p$ the $\ell_p$-norm for any $p\geq 1$. Given an integer $d >0$, the unit $\ell_p$-ball and $\ell_p$-sphere are respectively defined by
$$
\mathcal{B}_p := \left\{ \bo{x} \in \R^d \, : \, ||\bo{x}||_p  <1 \right\}; 
\qquad \quad 
\partial\mathcal{B}_p := \left\{ \bo{x} \in \R^d \, : \, ||\bo{x}||_p  = 1 \right\} \, .  
$$       
  
Given $i_k \in \N$ with $k=1, \ldots, d$, $\Vec{\boldsymbol{\imath}} :=(i_1, \ldots, i_d)$ and $\alpha \geq 1$, consider the cross-partial derivative operator $\mathcal{D}^{(\Vec{\boldsymbol{\imath}})}\M :=\prod_{k=1}^{d} \frac{\partial^{i_k} \M}{\partial^{i_k}  x_k}$, and define the H\"older  space of $\alpha$-smooth functions by         
$$         
\displaystyle              
\mathcal{H}_\alpha := \left\{    
\M : \R^d \to \R \, : \,
\left|\M(\bo{x}) - \sum_{0\leq i_1+\ldots+i_d \leq \alpha-1} 
\frac{\mathcal{D}^{(\Vec{\boldsymbol{\imath}})}\M(\bo{y})}{i_1! \ldots i_d!} \left(\bo{x}-\bo{y} \right)^{\Vec{\boldsymbol{\imath}}} \right| 
  \leq M_\alpha \norme{\bo{x}- \bo{y}}^\alpha \right\} \, ,        
$$            
with $M_\alpha>0$ a constant. Thus, $\nabla \M := \left[\frac{\partial \M}{\partial x_{1}}, \ldots, \frac{\partial \M}{\partial x_{d}} \right]^\T$ stands for the formal or Euclidean gradient of the $2$-smooth function $\M$, where $\frac{\partial \M}{\partial x_{k}}$ stands for the traditional partial derivative of $\M$ w.r.t. $x_k$, $k=1, \ldots, d$. \\ 
In what follows, $\esp[\cdot]$  and $\var[\cdot]$ stand for the expectation operator and the variance operator, respectively. 
   
\subsection{Non-Euclidean or dependent gradients} 
This section deals with the definition of the gradients of $2$-smooth functions evaluated at non-independent variables and the link between the dependent and the traditional or formal gradients.    \\

Denote with $\bo{\X} :=(X_1, \, \ldots,\, X_d)$ a random vector of continuous and non-independent variables having $F$ as the joint cumulative distribution function (CDF) (i.e., $\bo{\X}  \sim F$). For any $j \in \{1,\ldots, d\}$, $F_{x_j}$ or $F_j$ denotes the marginal CDF of $\X_j$, and  $F_j^{-1}$ stands for its  inverse. Also, denote $(\sim j) :=(1, \ldots, j-1, j+1,\ldots,  d)$ and $\bo{\X}_{\sim j} :=(\X_{1}, \ldots, \X_{j-1}, \X_{j+1}, \ldots, X_{d})$. The equality (in distribution) $\X \stackrel{d}{=} Z$ means that $\X$ and $Z$ have the same CDF.\\   
 
The formal $\nabla \M$ consists of the traditional partial derivative of $\M$ w.r.t. $x_k$'s, regardless the dependency structures of the input variables. For dependent variables $\bo{\X}$, there is a dependency function $r_j : \R^d \to \R^{d-1}$ such that (\cite{skorohod76,lamboni22b,lamboni21,lamboni25math,lamboni23mcap})   
$$
\bo{\X}_{\sim j}  = r_j\left(\X_j, \bo{Z}_{\sim j} \right) =: \left(r_{1,j}(\X_j, Z_1) , \ldots, r_{j-1,j}(\cdot),
 r_{j+1,j}(\cdot), \ldots, 
 r_{d,j}(\X_j, \bo{Z}_{\sim j})\right) \, ,          
$$
where $\bo{Z}_{\sim j} :=(Z_1, \ldots, Z_{j-1}, Z_{j+1}, \ldots, Z_d)$ is a random vector of independent variables, which is independent of $\X_j$. Such dependency functions are used for deriving the dependent partial derivatives of $\bo{x}$ w.r.t. $x_{j}$, that is, (\cite{lamboni21,lamboni23axioms})     
$$   
J^{(j)} \left(\bo{x} \right) := \frac{\partial \bo{x}}{\partial x_{j}} =       
 \left[\frac{\partial r_{1,j}}{\partial x_{j}}\, \ldots \, \underbrace{1}_{j^{\text{th}} \, \text{position}} \, \ldots \, \frac{\partial r_{d,j}}{\partial x_{j}} \right]^T \left(x_{j},\, r_{j}^{-1}\left(\bo{x}_{\sim j} |x_{j} \right) \right) \, ,                    
$$
and the dependent Jacobian matrix of $\bo{x} \mapsto \bo{x}$ given by (see \cite{lamboni21,lamboni23axioms} for more details)   
\begin{equation} \label{eq:grdxk}       
J^d \left(\bo{x} \right) := \left[J^{(1)} \left(\bo{x} \right), \ldots, J^{(d)} \left(\bo{x} \right) \right] \,  .    \nonumber                                  
\end{equation}      
Note that $J^d \left(\bo{x} \right)$ comes down to the identity matrix for independent variables. The dependent or non-Euclidean gradient of $\M$ with non-independent variables is given by (\cite{lamboni23axioms})     
\begin{equation} \label{eq:grad}
grad(\M)(\bo{x}) := G^{-1}(\bo{x})   \nabla\M(\bo{x})  \, ,    
\end{equation}     
with $ G(\bo{x}) :=  J^d \left(\bo{x} \right)^{\T} J^d \left(\bo{x} \right)$ the tensor metric and $G^{-1}(\bo{x})$ its generalized inverse. It is worth noting that the above derivation of gradients relies on the following assumption: $\bo{\X}$ consists of one block of dependent variables. In the presence of $K$ blocks of independent random vectors, the procedure is similar by treating independently such $K$ blocks, leading to the same form of the gradient (see \cite{lamboni23axioms} for more details). \\  

For independent variables, we have $G^{-1}(\bo{x}) =\mathcal{I}$ and $grad(\M)(\bo{x}) = \nabla\M(\bo{x})$. Since $G^{-1}(\bo{x})$ is always  known or can be computed, all of the results derived in this paper about the estimation or computation of $grad(\M)(\bo{x})$ are applicable to $\nabla\M(\bo{x})$ as well, and vice versa. For the sake of generality and precision, we are going to treat $grad(\M)(\bo{x})$, and the results for $\nabla\M(\bo{x})$ are immediate by taking $G^{-1}(\bo{x}) =\mathcal{I}$. 
   
\subsection{Stochastic surrogates of the dependent gradient and  estimators}
This section aims at extending the $L$-point-based surrogates of dependent gradients provided in \cite{lamboni24axioms} by considering the wide class of random vectors given by Equation (\ref{eq:proned}), including the $\ell_p$-spherical distributions. \\ 

Given $L, \, q \in \N\setminus\{0\}$, $\beta_{\ell} \in \R$ with $\ell=1, \ldots, L$, $\boldsymbol{\hh} := (\hh_1, \ldots, \hh_d) \in \R^d_+$, consider $L$ constraints given by
$ 
\sum_{\ell=1}^{L} C_{\ell} \beta_{\ell}^{r}= \delta_{1,r}
$
where  $r =0, \ldots, L-1$ or  $r =1, \ldots, L$   with $\delta_{1,r}$ the Kronecker symbol, that is, $\delta_{1,r}=1$ if $r=1$ and zero otherwise. One can see that the above  constraints lead to the existence of the constants $C_1, \ldots, C_L$ because some constraints rely on the Vandermonde matrix of the form   
$$             
A_{L} := 
\left[ \begin{array}{cccc}   1 & 1  & \ldots & 1\\  \beta_1 & \beta_2 & \ldots & \beta_{L} \\ 
 \beta_1^{2} & \beta_2^{2} & \ldots & \beta_{L}^{2} \\   
 \vdots &  \vdots  &  \vdots  & \vdots    \\ 
\end{array}  \right] \,  ,      
$$
which is invertible, as the determinant $\det\left(A_{L}\right) =\prod_{1\leq \ell_1 < \ell_2 \leq  L}\left( \beta_{\ell_1} - \beta_{\ell_2} \right)$ differs from zero for distinct $\beta_\ell$s (see \cite{lamboni24axioms,lamboni24stats} for more details). \\          
 
Denote with $\bo{W} :=(W_1, \ldots, W_d)$ a $d$-dimensional random vectors of independent variables satisfying: $\forall \, j \in \{1, \ldots, d\}$, 
$$   
 \esp\left[W_j\right] =0; \qquad \esp\left[\left(W_{j} \right)^{2}\right] =\sigma^2;   \qquad \esp\left[\left(W_{j} \right)^{2q+1}\right] =0; \qquad \esp\left[\left(W_{j} \right)^{2q}\right] < +\infty   \, .
$$    
It is shown in \cite{lamboni24axioms} (Theorem 1) that there exists  $\alpha_1 \in \{1, \ldots, L\}$ and reals coefficients $C_{1}, \ldots, C_{L}$ such that                 
\begin{equation}  \label{eq:approxful}           
\displaystyle     
grad(\M)(\bo{x}) =  G^{-1} (\bo{x})   \sum_{\ell=1}^{L} C_{\ell} \, \esp \left[ \M\left(\bo{x} + \beta_\ell \boldsymbol{\hh}\bo{W} \right) \frac{\bo{W} \bo{h}^{-1}}{\sigma^2} \right] +   \mathcal{O}\left( \norme{\bo{\hh}}^{2 \alpha_{1}} \right) \indic_{\bu} \, ,                   
\end{equation}            
with $\indic_{\bu} :=\left[1, \ldots, 1 \right]^T \in \R^d$, and 
 using the pointwise product $\boldsymbol{\hh}\bo{W}  := (\hh_1 W_1; \ldots, \hh_d W_d)$; $\boldsymbol{\hh}^{-1}\bo{W}  := (W_1/\hh_1; \ldots, W_d/\hh_d)$.  \\    

Corollary \ref{coro:stgradlp} extends such surrogates of $grad(\M)(\bo{x})$ by considering the wide class of random vectors $\bo{V} :=(V_1, \ldots, V_d)$ verifying   
\begin{equation}    \label{eq:proned}
\esp\left[V_k^2 \right]=\sigma^2,  \; \forall\, k \in \{1, \ldots, d \}; 
\qquad        
 \esp\left[\prod_{k=1}^d V_k^{q_k}\right] =0, \; \,\mbox{if, at least, one $q_k \in \N$ is odd} \, .  
\end{equation} 
Instances of random vectors that satisfy (\ref{eq:proned}) are listed below:     
\begin{itemize}
\item independent random variables such as $\bo{W}$;
\item random vectors that are uniformly distributed on the unit $p$-sphere (\cite{song97,valle12}); 
\item random vectors that are uniformly distributed over the unit $p$-ball (\cite{barthe05,barthe10,javid19}); 
\item random vectors that are $p$-spherically distributed on $\R^d$ (\cite{gupta97,valle12});
\item other new random vectors introduced in \cite{javid19,richter19}.    
\end{itemize}        
   
 Formally, denote with $\mathcal{S}_{d, p}$ a class of $p$-spherical distributions on $\R^d$, endowed with the $\ell_{p}$-norm when $p\geq 1$ (\cite{gupta97}) and the antinorm if $0< p <1$ (\cite{valle12}). The random vector $\bo{V} \sim \mathcal{S}_{d, p}$ is symmetrically distributed about zero.    
 Other properties of $\bo{V} \sim \mathcal{S}_{d, p}$ can be found in \cite{gupta97,valle12}, such as $\bo{V}$ admits a stochastic representation of the form  $\bo{V} =R \bo{U}$, where $R>0$ is a random variable, which is independent of  the generalized uniformly-distributed random vector on the unit ${p}$-sphere $\bo{U}$ (\cite{song97,valle12}). The $p$-generalized Gaussian distribution, and $p$-generalized Student distributions are examples. Remark that $\tilde{\bo{V}} := R \tilde{\bo{U}}$ also  satisfies (\ref{eq:proned}) with $\tilde{\bo{U}}$ a random vector that is uniformly distributed over the unit $p$-ball.              
\begin{corollary}  \label{coro:stgradlp} 
Let $\bo{V}$ satisfying (\ref{eq:proned}). Assume that $\M \in \mathcal{H}_\alpha$ with $\alpha \in \{2, \ldots, 2L+1\}$ and $\beta_\ell$s are distinct. Then, there exists  $\theta \in \{1, \ldots, L\}$ and reals coefficients $C_{1}, \ldots, C_{L}$ such that           
\begin{equation}  \label{eq:approxful}      
\displaystyle     
grad(\M)(\bo{x}) =  G^{-1} (\bo{x})   \sum_{\ell=1}^{L} C_{\ell} \, \esp \left[ \M\left(\bo{x} + \beta_\ell \boldsymbol{\hh}\bo{V} \right) \frac{\bo{V} \bo{h}^{-1}}{\sigma^2} \right] +   \mathcal{O}\left( \norme{\bo{\hh}}^{2 \theta} \right) \indic_{\bu} \, .                      
\end{equation}            
\end{corollary} 
\begin{preuve}
It is straightforward, as $\bo{W}$ and $\bo{V}$ share the same needed properties (i.e., (\ref{eq:proned})) for deriving the result. (see \cite{lamboni24axioms,lamboni24stats} and the supplementary documents for detailed proofs).   
\hfill $\square$    \\    
\end{preuve} 

The setting $L=1, \beta_1 =1, \, C_{1}=1$ or $\sum_{\ell=1}^{L=2} C_{\ell} \beta_{\ell}^{r}= \delta_{1,r}; \;  r =0, 1$
lead to the order of approximation $\mathcal{O}\left( \norme{\bo{\hh}}^{2} \right)$. For given distinct $\beta_\ell$s, the constraints 
$ 
\sum_{\ell=1}^{L} C_{\ell} \beta_{\ell}^{r}= \delta_{1,r}; \; r =1, 3, 5, \ldots, 2L-1$ allow for increasing that order up to $\mathcal{O}\left( \norme{\bo{\hh}}^{2L} \right)$.  
For the sake of simplicity, $\hh_j=\hh$ is used to define a neighborhood of a sample point of $\bo{X}$ (i.e., $\bo{x}$). Thus, using $\beta_{max} :=\max \left(|\beta_1|, \ldots, |\beta_L|\right)$ and keeping in mind the variance of $\beta_{\ell}\hh V_1$, it is reasonable to require $\beta_{max} \hh \sigma   \leq 1/2$. But, one can always choose $\sigma$ or $h$ to satisfy such a requirement (see Remark \ref{rem:chsig}). \\  

In view of Corollary \ref{coro:stgradlp},  the $L$-point-based surrogate of  $grad(\M)$ is given by
$$
\widetilde{grad(\M)}(\bo{x}) := G^{-1} (\bo{x})   \sum_{\ell=1}^{L} C_{\ell} \, \esp \left[ \M\left(\bo{x} + \beta_\ell \boldsymbol{\hh}\bo{V} \right) \frac{\bo{V} \bo{h}^{-1}}{\sigma^2} \right] \, , 
$$  
and one can compute the gradient using the method of moments.  Indeed, given an i.i.d. sample of $\bo{V}$, that is, $\left\{ \bo{V}_i :=(V_{i,1}, \ldots, V_{i,d})\right\}_{i=1}^N$, the consistent estimator of the dependent gradient is given by
\begin{equation}  \label{eq:estgrad}     
\displaystyle     
\widehat{grad(\M)}(\bo{x}) :=  \frac{G^{-1} (\bo{x})}{N \hh \sigma^2} \sum_{i=1}^N  \sum_{\ell=1}^{L} C_{\ell} \, \M\left(\bo{x} + \beta_\ell \boldsymbol{\hh}\bo{V}_i \right) \bo{V}_i  \, .         
\end{equation}       
Usually, $N=1$ for stochastic optimization problems. When $L=1$, $\widehat{grad(\M)}(\bo{x})$ is modified as follows:  
$$
 \frac{G^{-1} (\bo{x})}{N \hh \sigma^2} \sum_{i=1}^N  \left\{\M\left(\bo{x} + \boldsymbol{\hh}\bo{V}_i \right) - \overline{\M\left(\bo{x} + \boldsymbol{\hh}\bo{V} \right)}  \right\} \bo{V}_i \, , 
$$  
with  
$
\overline{\M\left(\bo{x} + \boldsymbol{\hh}\bo{V} \right)}  :=  \frac{1}{N} \sum_{i=1}^N  \M\left(\bo{x} + \boldsymbol{\hh}\bo{V}_i \right)    
$.     
         
\begin{rem}
The estimator $\widehat{grad(\M)}(\bo{x})$ is an unbiased estimator of the $L$-point-based surrogate of  $grad(\M) 
$, that is, 
$ \esp\left[\widehat{grad(\M)}(\bo{x}) \right] = \widetilde{grad(\M)}(\bo{x}) 
$. \\
For one-point stochastic optimization problems,  $\overline{\M\left(\bo{x} + \boldsymbol{\hh}\bo{V} \right)}$ may be considered null at the beginning and updated at each step $t$ as follows:  
$
\overline{\M\left(\bo{x} + \boldsymbol{\hh}\bo{V} \right)}  :=  \frac{1}{t} \sum_{i=1}^t  \M\left(\bo{x} + \boldsymbol{\hh}\bo{V}_t \right)   
$. 
\end{rem}  
 
In what follows, we will rely on $\bo{V} \sim \mathcal{S}_{d, p}$  or $\tilde{\bo{V}} := R \tilde{\bo{U}}$ to derive the main results. Recall that $\tilde{\bo{U}}$ is a random vector that is uniformly distributed over the unit $p$-ball, and there is a bijection between $\bo{U}$ and $\tilde{\bo{U}}$. Moreover, $\bo{V} \sim \mathcal{S}_{d, p}$ and $\tilde{\bo{V}} := R \tilde{\bo{U}}$ share the same exponential concentration inequalities (\cite{ledoux01,louart20}).           
            
\section{Convergence analysis} \label{sec:dfexpderi}
This section deals with the statistical properties of the estimator of the gradient given by Equation (\ref{eq:estgrad}), such as its mean squared error (MSE) and rate of convergence. 
\subsection{Dimension-free upper-bounds of the bias} \label{sec:dfparder} 
It is worth mentioning that dimension-free upper-bounds of the bias have been provided in \cite{lamboni24axioms} under different structural assumptions on the  deterministic functions $\M$ using $\bo{W}$. This section aims at deriving the dimension-free upper-bounds of the bias using $\bo{V} \sim \mathcal{S}_{d,p}$ and the minimal requirement about the smoothness of deterministic functions. 
\begin{assump}[A1]
 The deterministic function $\M \in \mathcal{H}_{\alpha}$ for any $\alpha \in \{1, 2\}$.   
\end{assump}
 Corollary~\ref{coro:parderord} gives the upper-bounds of the bias when $L=2$ and $\hh_j=\hh$. For a given matrix $\mathcal{M}$, denote with $\left|\mathcal{M}\right|$ the matrix whose entries are the  absolute values of those of $\mathcal{M}$; and define   
$$
K_{1,d,p} := 
\frac{\Gamma(d/p) \,\left[\Gamma(4/p) \Gamma(1/p) + (d-1)\Gamma(3/p) \Gamma(2/p)\right]}{\Gamma^2(1/p) \Gamma(d/p+3/p) } \, ,  
$$    
with $\Gamma(\cdot)$ the Gamma function.               
       
\begin{corollary} \label{coro:parderord}         
Let $L=2$, $\beta_1 =-\beta_2 =1$ and $C_{1} =-C_{2}= 1/2$. If (A1) holds, then there  exists $M_2 >0$ such that      
\begin{equation}  \label{eq:bias}                          
\displaystyle          
\norml{
 grad(\M)(\bo{x}) -  G^{-1} (\bo{x})   \sum_{\ell=1}^{L=2} C_{\ell} \, \esp \left[ \M\left(\bo{x} + \beta_\ell \hh \bo{V} \right) \frac{\bo{V}}{\sigma^2 \hh} \right]} \leq  \frac{ M_2 \hh \, \esp\left[ R^3 \right]}{\sigma^2} \norml{ \left|G^{-1}(\bo{x})\right| \indic_{\bu}} K_{1,d,p}  \, ;              
\end{equation}     
\begin{equation}  \label{eq:biasl2}                                
\displaystyle          
  \norme{
 grad(\M)(\bo{x}) -  G^{-1} (\bo{x})   \sum_{\ell=1}^{L=2} C_{\ell} \, \esp \left[ \M\left(\bo{x} + \beta_\ell \hh \bo{V} \right) \frac{\bo{V}}{\sigma^2 \hh} \right]} \leq  \frac{ M_2 \hh \, \esp\left[ R^3 \right]}{\sigma^2} \norme{ \left|G^{-1}(\bo{x})\right| \indic_{\bu}} K_{1,d,p}    \, .     
\end{equation} 
\end{corollary}          
\begin{proof}    
See Appendix \ref{app:coro:parderord}.  
\end{proof} 

\begin{rem}  \label{rem:l1}
The results derived in Corollary \ref{coro:parderord} remain valid for the setting $L=1$, $C_1=\beta_1=1$ (see Appendix \ref{app:coro:parderord}).   
\end{rem}

Note that the above results hold for any random variable $R>0$. Since $\esp\left[V_k^2\right] =\esp\left[V_1^2\right] = \esp\left[U_1^2\right] \esp\left[R^2\right] =\sigma^2$, one can deduce     
$
\esp\left[R^2\right] = \sigma^2 \frac{\Gamma(1/p) \Gamma(d/p+2/p)}{\Gamma(3/p) \Gamma(d/p)} 
$    
(see \cite{valle12}, Corollary 2.8).  An interesting choice of $R$ will lead to dimension-free upper-bound of the bais for any $p>0$, and a small error of approximations. Thus, taking $R_0 \sim \mathcal{U}(0,\, \xi)$ requires setting $\xi := \sqrt{3 \esp\left[R^2\right]}$. To be able to derive dimension-free mean squared errors, it is essential to exhibit the dimensionality in the expression of $\esp\left[R_0^q\right]$ for any $q \in \N$.   

\begin{lemma}  \label{lam:qthm}
Let $R_0 \sim \mathcal{U}(0,\, \xi)$ with $\xi := \sqrt{3 \esp\left[R^2\right]}$. \\

$\quad$ (i)   If $1 \leq p \ll d$, then 
\begin{equation} \label{eq:rmo} 
\esp\left[R_0^q\right] =  \frac{3^{q/2} \sigma^q}{q+1}
\left( \frac{\Gamma(1/p) \Gamma(d/p+2/p)}{\Gamma(3/p) \Gamma(d/p)}  \right)^{q/2}
\approx \frac{3^{q/2} \sigma^q \, d^{q/p}}{(q+1) \, p^{q/p}}
\left( \frac{\Gamma(1/p)}{\Gamma(3/p)} \right)^{q/2}
 \, ,   
\end{equation}   
    
$\quad$ (ii)   If $1 \leq d \ll p$, then 
\begin{equation} \label{eq:rmo1} 
\esp\left[R_0^q\right] =  \frac{3^{q/2} \sigma^q}{q+1}
\left( \frac{\Gamma(1/p) \Gamma(d/p+2/p)}{\Gamma(3/p) \Gamma(d/p)}  \right)^{q/2}
\approx  \frac{3^{q} \sigma^q}{q+1}
\left( \frac{d}{d+2} \right)^{q/2}  \, .         
\end{equation}     
\end{lemma}
\begin{proof}
See Appendix \ref{app:lam:qthm}.      
\end{proof}
 
In view of Equation (\ref{eq:rmo1}), the $q$th-order moment of $R_0$ does not depend on $d$ in higher-dimensions. Indeed, when $d, p \to +\infty$ with $d/p \to 0$, we have $\esp\left[R_0^q\right]
\approx  \frac{3^{q} \sigma^q}{q+1}$. Also, one can see that $\esp\left[R_0^q\right] \lessapprox  \frac{3^{q} \sigma^q}{q+1}$ in general, where $\lessapprox$ stands for the less approximation, leading to the approximated upper-bound. Based on $R_0$ and its $q$th-order moment, the dimension-free bias is derived in the following corollary.  For that purpose,  consider   
$$
K_{2,d,p}  :=  \frac{3\sqrt{3}}{4}
\frac{\left[\Gamma(4/p) \Gamma(1/p) + (d-1)\Gamma(3/p) \Gamma(2/p)\right] \Gamma^{3/2}(d/p+2/p)}{ \Gamma^{1/2}(d/p) \Gamma^{1/2}(1/p) \Gamma(d/p+3/p)  \Gamma^{3/2}(3/p)} \, ,  
$$
which is approximated  by    
$$
K_{2,d,p}  \approx  \left\{
\begin{array}{cl} 
  \frac{3\sqrt{3}}{4}  
\frac{\left[\Gamma(4/p) \Gamma(1/p) + (d-1)\Gamma(3/p) \Gamma(2/p)\right]}{\Gamma^{1/2}(1/p) \Gamma^{3/2}(3/p)} &  \mbox{if}\, 1\leq p \ll d \\
\frac{9(d+3)(2d+1) d^{1/2}}{16(d+2)^{3/2}}  &  \mbox{if}\, 1\leq d \ll p \\
\end{array}
\right. \,  .                
$$

\begin{corollary} \label{coro:dimfreup}         
Let $L=2$, $\beta_1 =-\beta_2 =1$ and $C_{1} =-C_{2}= 1/2$. If (A1) holds and $R =R_0$,  then       
\begin{equation}                                  
\displaystyle          
  \norme{
 grad(\M)(\bo{x}) -  G^{-1} (\bo{x})   \sum_{\ell=1}^{L=2} C_{\ell} \, \esp \left[ \M\left(\bo{x} + \beta_\ell \hh \bo{V} \right) \frac{\bo{V}}{\sigma^2 \hh} \right]} \leq  M_2 \hh \sigma \norme{ \left|G^{-1}(\bo{x})\right| \indic_{\bu}} K_{2,d,p}    \, .        \nonumber
\end{equation} 
Moreover, if $\sigma =\left( \norme{ \left|G^{-1}(\bo{x})\right| \indic_{\bu}} K_{2,d,p} \right)^{-1}$, then
\begin{equation}  \label{eq:biasl2df}                                    
\displaystyle          
  \norme{
 grad(\M)(\bo{x}) -  G^{-1} (\bo{x})   \sum_{\ell=1}^{L=2} C_{\ell} \, \esp \left[ \M\left(\bo{x} + \beta_\ell \hh \bo{V} \right) \frac{\bo{V}}{\sigma^2 \hh} \right]} \leq  M_2 \hh     \, .        
\end{equation}       
\end{corollary}           
\begin{proof}    

See Appendix \ref{app:coro:dimfreup}.    
\end{proof}  

It appears that the proposed approach leads to the upper-bound of the bias  that does not depend on the dimensionality for any value of $p>0$. We also have  
\begin{equation}  \label{eq:biasdf}                          
\displaystyle           
\norml{
 grad(\M)(\bo{x}) -  G^{-1} (\bo{x})   \sum_{\ell=1}^{L=2} C_{\ell} \, \esp \left[ \M\left(\bo{x} + \beta_\ell \hh \bo{V} \right) \frac{\bo{V}}{\sigma^2 \hh} \right]} \leq  M_2 \hh   \, ;              
\end{equation}   
by taking $\sigma = \left(\norml{ \left|G^{-1}(\bo{x})\right| \indic_{\bu}} K_{2,d,p} \right)^{-1}$. 
    
\begin{rem} \textbf{Dimension-free upper-bounds for $L$-point-based surrogates}. \\
Note that the above dimension-free upper-bound (i.e., $M_2 \hh$) has been already obtained in \cite{lamboni24axioms} using $\bo{W}$. Such a bound still holds for any radial variable $R$, provided that $\esp\left[R^3\right]$ is proportional to $\sigma^3$ (see Corollary~\ref{coro:parderord}). Of course, a trivial choice is given by $R_1 \sim \delta_{\sigma}$ with $\delta_{\sigma}$ the Dirac probability measure. \\ 
One can check that dimension-free bias holds using $L=1$ (see Remark \ref{rem:l1}). For the sequel of generality, in the presence of highly smooth functions, increasing $L$ or equivalently the number of evaluations of $\M$ at randomized points will end up with the same dimension-free upper-bound, except the choice of the value of $\sigma$ (see \cite{lamboni24axioms} for more details).  
\end{rem}     
                
\begin{rem} \textbf{Choice of $\sigma$ for any $L \in\{1, 2 \}$.} \label{rem:chsig} \\
In the case of independent inputs, one can see that taking $\sigma \in\{ d^{3/2}, \; d^{-2}\}$ is sufficient to obtain dimension-free upper-bounds of the bias. This choice remains valid for dependent variables because $\norme{ \left|G^{-1}(\bo{x})\right| \indic_{\bu}}\leq  \sqrt{d} \, ||G^{-1}(\bo{x}) ||_s$ with $||\cdot ||_s$ the spectral norm.   
\end{rem}     
   
\subsection{Mean squared errors of estimators of dependent gradients} \label{sec:dfderest}  
Using a sample of $\bo{V} \sim \mathcal{S}_{d,p}$, that is, 
$\left\{ \bo{V}_i \right\}_{i=1}^N$, recall that the estimator of $grad(\M)(\bo{x})$ is  given by 
$$ 
\displaystyle    
\widehat{grad(\M)}(\bo{x}) = \left\{ \begin{array}{cc} 
 \frac{G^{-1}(\bo{x})}{N \sigma^2 \hh}  \sum_{i=1}^ N \sum_{\ell=1}^{L} C_{\ell}  \M\left(\bo{x} + \beta_\ell \hh\bo{V}_i \right) \bo{V}_i  & \mbox{if L>1} \\
 \frac{G^{-1}(\bo{x})}{N \sigma^2 \hh}  \sum_{i=1}^ N \left\{ \M\left(\bo{x} + \hh\bo{V}_i \right)- \overline{\M\left(\bo{x} + \hh\bo{V} \right)} \right\} \bo{V}_i & \mbox{if L=1} \\
\end{array} 
\right.  \, .  
$$ 
 
Before providing the MSEs, the following intermediate results are needed. 
Recall that $\bo{V} = R \bo{U}$, where $\bo{U}$ is the generalized uniformly-distributed random vector on the unit $p$-sphere (\cite{song97}). It is also shown in \cite{ledoux01,louart20} that $\bo{U}$ is a $p$-exponential concentrated random vector as well as $g(\bo{U})$ for any Lipschitz function $g$. Based on these elements, the $q$th-order moments of $g(\bo{V})$ are given below.        
  
\begin{lemma} \label{lem:mosp}
Let $\bo{V} \sim \mathcal{S}_{d, p}$  be $\ell_{p}$-spherically distributed with $p\geq 1$ or $\bo{V} \stackrel{d}{=}\tilde{\bo{V}}$ and $g$ be a $L_0$-Lipschitz function w.r.t. $\bo{V}$. If $g(\bo{V})$ has finite $q$th-order moments, then      
$$
\esp\left[\left|g\left(\bo{V} \right)- g\left(\bo{0} \right)\right|^q\right] \leq
\frac{q C  \Gamma(q/p)}{p c_0^{q/p}} \frac{L_0^q}{d^{^{q/p}}} \esp\left[R^{q} \right]  \, ,        
$$    
where $C, c_0$ are constants depending only on $p$.  
Moreover,     
$$ 
\esp\left[ R^q  \left|g\left(\bo{V} \right)- g\left(\bo{0} \right)\right|^q\right] \leq
\frac{q C  \Gamma(q/p)}{p c_0^{q/p}} \frac{L_0^q}{d^{^{q/p}}} \esp\left[R^{2q} \right]  \, .          
$$              
\end{lemma}             
\begin{proof} 
See Appendix \ref{app:lem:mosp}.     
\end{proof}
  
For particular choices of the radial variable $R$, the following precise results in terms of the power of the dimensionality are obtained.       
\begin{corollary} \label{coro:mosp}  
Under the conditions of Lemma \ref{lem:mosp},  consider  $\bo{V} \sim \mathcal{S}_{d, p}$ or $\tilde{\bo{V}}$  with $R \stackrel{d}{=} R_0$. \\

$\quad$ (i) If $1\leq p \ll d$, then        
$$
\esp\left[\left|g\left(\bo{V} \right)- g\left(\bo{0} \right)\right|^q\right]  \lessapprox
\frac{3^{q/2} \sigma^q  q C  \Gamma(q/p) L_0^q}{(q+1) p^{1+q/p} c_0^{q/p} } \left( \frac{\Gamma(1/p)}{\Gamma(3/p)} \right)^{q/2}  \, ;    
$$    
$$
\esp\left[ R^q  \left|g\left(\bo{V} \right)- g\left(\bo{0} \right)\right|^q\right]  \lessapprox 
\frac{3^q \sigma^{2q} q C  \Gamma(q/p) L_0^q \, d^{^{q/p}}}{(2q+1) p^{1+2q/p} c_0^{q/p} } \left( \frac{\Gamma(1/p)}{\Gamma(3/p)} \right)^{q} \, .                      
$$ 

$\quad$ (ii) If $1\leq d \ll p$ and $q \ll p$, then            
$$
\esp\left[\left|g\left(\bo{V} \right)- g\left(\bo{0} \right)\right|^q\right]  \lessapprox
\frac{3^{q} \sigma^q C L_0^q}{(q+1) c_0^{q/p}}\left( \frac{d}{d+2} \right)^{q/2}  d^{-q/p}   \, .         
$$   
\end{corollary}   
\begin{proof}
See Appendix \ref{app:coro:mosp}. 
\end{proof}
        
For uniformly-distributed random vectors on the unit ${p}$-sphere, the following lemma helps for controlling the $q$th-order moment of the Euclidean norm of such variables.  
					   
\begin{lemma} \label{lem:molpuni}  
Let $\bo{U}$  be uniformly-distributed on the unit $p$-sphere or ball; $C'_0, c_0, c_0'$ be constants depending on $p$ and $p\geq 1$. Then, 
 
$$
\esp\left[ \left|\left|  \bo{U} \right|\right|_2^q \right]  \leq  \frac{qC'_0\Gamma(q/p)}{p(c_0')^{q/p}}; 
\qquad \quad 
\esp\left[ \left|\left|  \bo{U} \right|\right|_2^4 \right]  \leq  \frac{p^{4/p}\Gamma(5/p)}{\Gamma(1/p)} d^{2-4/p}
\, .     
$$             
\end{lemma} 			
\begin{proof}
See Appendix \ref{app:lem:molpuni}.  
\end{proof}		  
                     
Now, we have all the elements in hand to derive the MSEs of the gradient estimators and the corresponding rates of convergence. Theorem \ref{theo:mse} and Corollary \ref{coro:optderi1} provide such results. To that end, define 
$$
K_{3,p} :=  \frac{18 M_1^2}{5}  \frac{C^{1/2}  \Gamma^{1/2}(4/p) \Gamma^{1/2}(5/p) \Gamma^{3/2}(1/p)}{c_0^{2/p}  p^{1/2+2/p} \Gamma^2(3/p)}   
\left|\left| G^{-1}(\bo{x}) \right|\right|_s^2 \, ,   
$$   
$$
K_{4,p} := \frac{36 M_1^2}{5}  \frac{(C')^{1/2}  \Gamma(4/p) \Gamma^2(1/p)}{(c')^{2/p}  p^{1+4/p} \Gamma^2(3/p)} 
\left|\left| G^{-1}(\bo{x}) \right|\right|_s^2  \, ,
$$
with $||\cdot||_s$ the spectral norm of matrices. For independent variables, note that $\left|\left| G^{-1}(\bo{x}) \right|\right|_s =1$.         

\begin{assump}[A2]
The randomized function $\M(\bo{x} + \beta_\ell \hh \bo{V})$ has finite fourth-order moment.     
\end{assump}     
           
\begin{theorem}      \label{theo:mse} 
Let $L=2$, $\beta_1 =-\beta_2 =1$; $C_{1} =-C_{2}= 1/2$ and  $\bo{V} \sim \mathcal{S}_{d, p}$ with $R \stackrel{d}{=} R_0$. Assume $\sigma =\left( \norme{ \left|G^{-1}(\bo{x})\right| \indic_{\bu}} K_{2,d,p} \right)^{-1}$ and (A1)-(A2) hold. \\    

 $\quad$ (i) If $1 \leq p \ll d$,  then        
\begin{equation}  \label{eq:mse}                             
 \esp \left[\norme{ \widehat{grad(\M)}(\bo{x}) - grad(\M)(\bo{x})}^2 \right] \leq   M_2^2 \hh^2 +  
 N^{-1} \min\left( K_{3,p}\, d  ,\;  K_{4,p}\, d^{2/p} \right)   \, .           
\end{equation}          
      
 $\quad$ (ii) Let $K_{5,p} := \frac{81 C' M_1^2}{5} $. If $1 \leq d \ll p < +\infty$,  then        
\begin{equation}  \label{eq:mse2}                             
 \esp \left[\norme{ \widehat{grad(\M)}(\bo{x}) - grad(\M)(\bo{x})}^2 \right] \leq   M_2^2 \hh^2 +  
 N^{-1} K_{5,p}\, \frac{d^{2-2/p}}{(d+2)^2}  \, .           
\end{equation}                                   
\end{theorem}                                    
\begin{proof}     
See Appendix \ref{app:theo:mse}.            
\end{proof}    
 
Remark that the second terms of the upper-bound given by (\ref{eq:mse}) is linear with the dimensionality $d$ for every $p\geq 1$, which is the worst case, because $d^{2/p} \leq d$ for any $p$ verifying $2 \leq p \ll d$. When $d \ll p < +\infty$, one obtains $d^{-2/p} \approx 1$ in higher-dimensions (see Equation (\ref{eq:mse2})).  
Also, such second terms do not depend on the bandwidth $\hh$ and $\sigma$. This key observation leads to the derivation of the optimal and parametric rates of convergence of the estimator of the gradient. 
    
\begin{corollary}      \label{coro:optderi1} 
Under the conditions of Theorem \ref{theo:mse}, let $\hh \propto N^{-\gamma/2}$ with $\gamma \in ]1, \, 2[$. \\
  
$\qquad$ (i) If $1\leq p\leq 2$, then   $\esp \left[\norme{ \widehat{grad(\M)}(\bo{x}) - grad(\M)(\bo{x})}^2 \right] = \mathcal{O}\left(N^{-1} d \right)$.        
          
$\qquad$ (ii) If $2\leq p \ll d$, then $\esp \left[\norme{ \widehat{grad(\M)}(\bo{x}) - grad(\M)(\bo{x})}^2 \right] = \mathcal{O}\left(N^{-1} d^{2/p} \right)$.   
        
$\qquad$ (iii) If $d \ll p <+\infty$, then $\esp \left[\norme{ \widehat{grad(\M)}(\bo{x}) - grad(\M)(\bo{x})}^2 \right] = \mathcal{O}\left(N^{-1} d^{-2/p} \right)$,  
provided that $K_{5,p} \ll +\infty$.                   
\end{corollary}                           
\begin{proof}     
It is straightforward since $\hh^2 \propto N^{-\gamma}$ and $N \hh \to \infty$ when $N \to~\infty$. 
\end{proof} 

\begin{rem}  \textbf{Choice of $\hh$ for any $L \in\{1, 2 \}$.} \\ 
Choosing $\hh \propto N^{-\gamma/2}$ with $\gamma \in ]1, \, 2[$ is a consequence of i) the requirement $N \hh \to \infty$ when $N \to~\infty$ often used in non-parametric estimations, and ii) the willing to achieve the parametric rate of convergence. But, on can choose $\hh \leq N^{-1}$ to reduce the bias without modifying the parametric rate of convergence.      
\end{rem}

The derived rate of convergence depends on $d^{2/p}$, meaning that our estimator suffers from the curse of dimensionality for some values of $p$. Additionally to the dimension-free bias obtained, it turns out that taking $2< p \ll d$ allows for breaking down the curse of dimensionality (see Corollary \ref{coro:optderi2}), and for improving the best  known rates of convergence  (see \cite{bach16,shamir17,lamboni24axioms,berahas22,akhavan24,gasnikov23}). 
In higher-dimensions, taking higher values of $p$ could help for reaching dimension-free upper-bound of the MSE of the proposed estimator (see Point (iii) of Corollary \ref{coro:optderi1}). Indeed, for the setting $1\leq d \ll p<+\infty$, the rate of convergence obtained is  $N^{-1} d^{-2/p} \approx N^{-1}/e^2$, provided that the corresponding $K_{5,p}$ is not too large.   
  
\begin{corollary}      \label{coro:optderi2}  
Under the conditions of Theorem \ref{theo:mse}, let $\hh \propto N^{-\gamma/2}$ with $\gamma \in ]1, \, 2[$.\\
  If $ max\left(2, \log(d) \right) < p \ll d$ and $K_{4,p} \ll +\infty$, then        
$$
\esp \left[\norme{ \widehat{grad(\M)}(\bo{x}) - grad(\M)(\bo{x})}^2 \right] = \mathcal{O}\left(N^{-1} \right) \, . 
$$                 
\end{corollary}                               
\begin{proof}     
One can see that $d^{2/p} = \exp\left( \frac{2\log(d)}{p}\right) \leq \exp(2)$. 
\end{proof} 
It turns out that the dimension-free MSE of the proposed estimator is reached, provided that $K_{4,p_0}$ is not too large for the selected $p=p_0$. Since taking $[\max(2, log(d)]+1 \leq p \ll d$ help obtaining the optimal, parametric and dimension-free MSE, such values should be used in practice. Note that 
  $[\max(2, log(d)]$ stands for the largest integer that is less than  $\max(2, log(d)$. For instance, if $d=100$, we may use $p=5; 6$, and if $d=1000$, then $p=7; 8$.      
   
\begin{rem} 
The above results can be established using $\tilde{\bo{V}}$ as well thanks to Corollary \ref{coro:mosp} and Lemma \ref{lem:molpuni}, but with slight modifications of the proofs . For computational issues (see Section \ref{sec:fdm}), $\tilde{\bo{V}}$ offers more flexibility when applying the Gram-schmidt procedure.      
\end{rem}    
                                   
\section{Computational issues}\label{sec:fdm}  
This section provides simulated results using the proposed approach compared to well-established methods, such as i) the FDM using the R-package numDeriv (\cite{gilbert19}) with $\hh=10^{-4}$, ii) the MC approach provided in \cite{patelli10} with $\hh=10^{-4}$, and iii) the MC approach based on independent uniform distributions $\mathcal{U}(-1, 1)$ proposed in \cite{lamboni24axioms}. \\  

 The R-package LHS and the R-package greybox are used for generating the values of the $d$ independent $p$-generalized (standard) Gaussian variables (i.e., $\bo{G} :=(G_1, \ldots, G_d$)). Such values are then used to form $V_j =G_j/||\bo{G}||_p; \, j=1 \ldots, d$. Finally, as $\esp\left[V_{j_1} V_{j_2}\right]$ is estimated by $\frac{1}{N}\sum_{i=1}^N V_{i,j_1} V_{i,j_2}$ for a given sample size $N$, the Gram-schmidt procedure is applied to $\bo{V}$ (when possible) to ensure that $\esp\left[V_{j_1} V_{j_2}\right]=0$ for all $j_1, j_2 \in \{1, \ldots, d\}$ and $j_1 \neq j_2$.  
For higher values of $p$ (i.e., $p > 2000$), the corresponding $p$-generalized (standard) Gaussian distribution is approximated by the uniform distribution $\mathcal{U}(-1, 1)$ like in \cite{lamboni24axioms}.  \\    
To assess the numerical accuracy of each approach, the following error measure is considered:  
$$
Err := \frac{\norme{G^{-1}(\bo{0})\left[\nabla \M(\bo{0})- \widehat{\nabla \M}(\bo{0}) \right]}}{\norme{G^{-1}(\bo{0})\nabla \M(\bo{0})}} \, ,   
$$          
with $\widehat{\nabla \M}(\bo{0})$ the estimated value of the traditional gradient. In this section,  $\hh=10^{-4}$ and $\sigma=1/d^2$ are used.  

\subsection{Rosenbrock's function}
The Rosenbrock function is defined as follows: $\forall\, \bo{x} \in \R^d$,   
$$
r(\bo{x}) := \sum_{k=1}^{d-1} \left[ (1-x_k)^2 +  100 \left(x_{k+1} -x_k^2 \right)^2 \right] \, ,
$$     
and its traditional gradient at $\bo{0}$ is $\nabla r(\bo{0}) =\left[-2, \ldots, -2, 0 \right]^\T \in \R^{100}$ (see \cite{patelli10}). The dimensions  $d \in\{10, 100, 1000\}$ are considered. In the case of $d=10$,  independent input variables and correlated inputs are considered.  The correlation between $\X_{j_1}$ and $\X_{j_2}$  is fixed at $\mathcal{R}_{j_1 j_2} := \left(\frac{1}{2} \right)^{|j_1-j_2|}$ and the variance $\var\left[ \X_{j}\right] =1, j=1, \ldots, d$. Tables \ref{tab:fdm1}-\ref{tab:fdm4} report the values of $Err$ for the Rosenbrock function and for different approaches. 

\begin{table}[htbp]  
\begin{center}   
\begin{tabular}{lcccc}   
\hline      
\hline          
 \textbf{ d=10} $G =\mathcal{R}\mathcal{R}$   &     \multicolumn{4}{c}{Number of total model evaluations (i.e., $LN$)} \\   
				&   $11$   &   $15$  &  $20$ & $20$  \\        			
 Methods &      &     &   &  \\      
\hline 
				 &   $L=1$   &   $L=1$  &   $L=1$  &  $L=2$   \\   
  MC (\cite{lamboni24axioms}), $h=10^{-4}$  &   0.091    & 0.067     &  0.05   & 0.09  \\ 
 	This paper ($p=3$)  &   0.089     &  0.066     &  0.05  &     0.091  \\         
\hline                                   
\hline                       
    \end{tabular}                 
  \end{center}    
  \caption{Values of $Err$ for different approximations of $grad r(\bo{0})$ when $d=10$.}
			 	\label{tab:fdm1}      
\end{table}  
 
\begin{table}[htbp]  
\begin{center}   
\begin{tabular}{lcccc}   
\hline      
\hline          
 \textbf{ d=10} $G=\mathcal{I}$   &     \multicolumn{4}{c}{Number of total model evaluations (i.e., $LN$)} \\   
				&   $11$   &   $15$  &  $20$ & $20$  \\      			
 Methods &      &     &   &  \\      
\hline	
 FDM  (\cite{gilbert19})  &   -  &  -   &  - &   0.005   \\
  MC (\cite{patelli10}) &  0.0083  &  -    &  -    &   -  \\ 
\hline 
				 &   $L=1$   &   $L=1$  &   $L=1$  &  $L=2$   \\   
	This paper ($p=3$)  &   0.091     & 0.067     &  0.05  &     0.091  \\         
\hline                                         
\hline                          
    \end{tabular}                
  \end{center}    
  \caption{Values of $Err$ for different approximations of $\nabla r(\bo{0})$ when $d=10$.}
			 	\label{tab:fdm2}  
\end{table}        
      
\begin{table}[htbp]  
\begin{center}   
\begin{tabular}{lcccc}   
\hline      
\hline          
 \textbf{ d=100} $G=\mathcal{I}$   &     \multicolumn{4}{c}{Number of total model evaluations (i.e., $LN$)} \\   
				&   $101$   &   $150$  &  $200$ & $200$  \\      			
 Methods &      &     &   &  \\      
\hline	
 FDM  (\cite{gilbert19})  &   -  &  -   &  - &   0.005   \\
  MC (\cite{patelli10}) &    $0.042$  &  -    &  -    &   -  \\ 
\hline 
				 &   $L=1$   &   $L=1$  &   $L=1$  &  $L=2$   \\   
  MC (\cite{lamboni24axioms})  &   0.035     &  0.014     &  0.009  &     0.009  \\ 
		  MC (\cite{lamboni24axioms}), $h=10^{-4}$ &   0.0099     &  0.0067     &  0.005  &     0.009  \\ 
	  This paper ($p=5$)  &   0.0099     &  0.0066     &  0.005  &     0.0099  \\         
\hline                                   
\hline                            
    \end{tabular}                 
  \end{center}    
  \caption{Values of $Err$ for different approximations of $\nabla r(\bo{0})$ when $d=100$.}
			 	\label{tab:fdm3}        
\end{table}    

\begin{table}[htbp]  
\begin{center}   
\begin{tabular}{lcccc}   
\hline      
\hline          
 \textbf{ d=1000} $G=\mathcal{I}$   &     \multicolumn{3}{c}{Number of total model evaluations (i.e., $LN$)} \\   
				&   $1001$   &     $2000$ & $2000$  \\      			
 Methods &      &    &  \\      
\hline	
 FDM  (\cite{gilbert19})  &   -    &  - &   0.0051 \\
  MC (\cite{patelli10}) &    $0.036$     &  -    &   -  \\ 
\hline 
				 &   $L=1$    &   $L=1$  &  $L=2$   \\   
  MC (\cite{lamboni24axioms}), $h=10^{-4}$  & 0.0009   &  0.0005   & 0.0009 \\ 
	  This paper ($p=7$)  &  0.0015  &  0.0005   &  0.0015  \\                      
\hline                               
\hline                       
    \end{tabular}                   
  \end{center}    
  \caption{Values of $Err$ for different approximations of $\nabla r(\bo{0})$ when $d=1000$.}
			 	\label{tab:fdm4}           
\end{table}          
         
\begin{figure}[!hbp]            
\begin{center}
\includegraphics[height=13cm,width=9cm,angle=270]{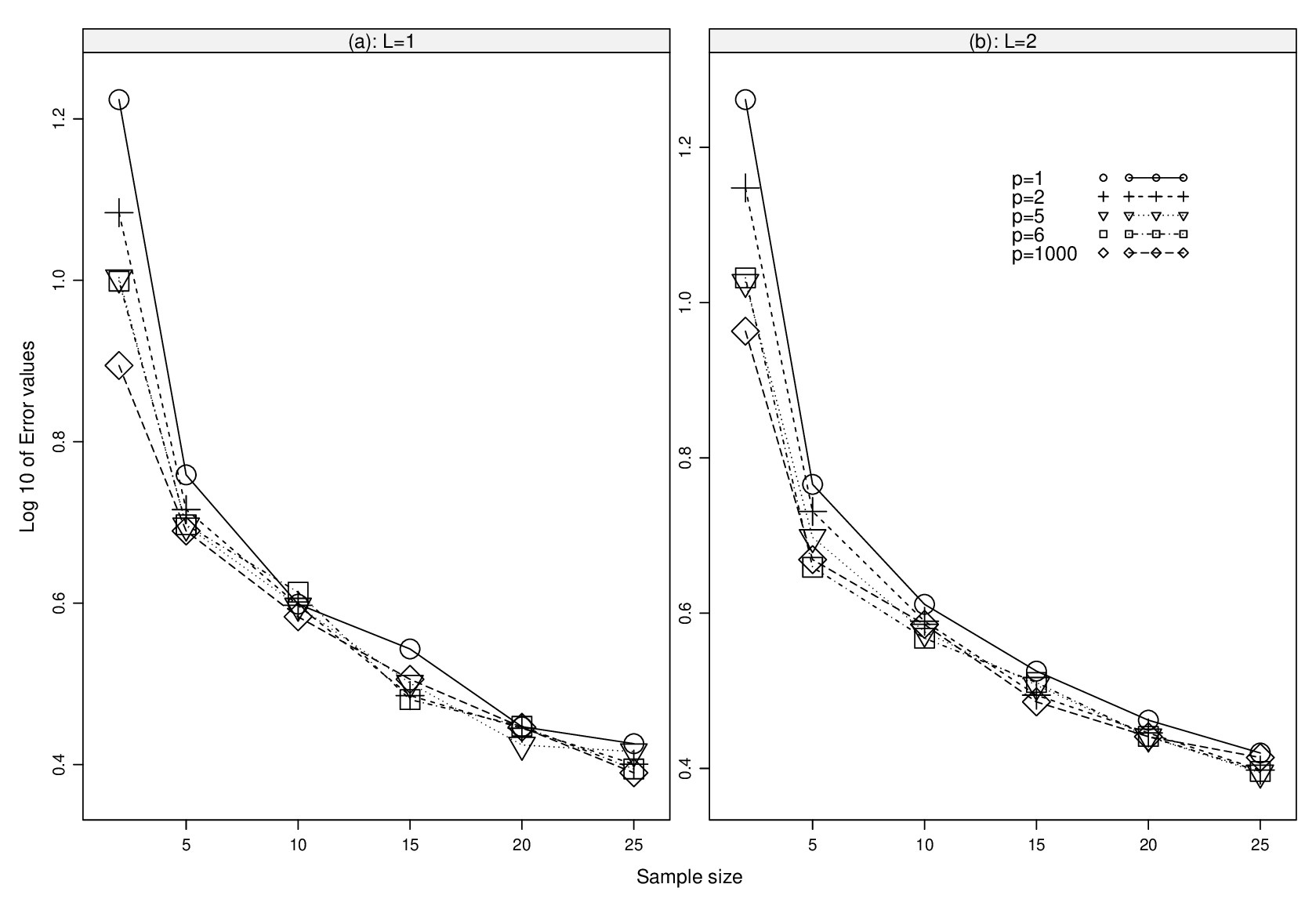}    
\end{center} 
\caption{Average of error values against the sample sizes ($N$) for $p$-spherical distributions (i.e., $\bo{V}$) when $\bo{d=100}$. Panel (a) is obtained when $L=1$, while Panel (b) is associated with $L=2$.} 
 \label{fig:gradsp100}         
\end{figure}

\begin{figure}[!hbp]         
\begin{center}
\includegraphics[height=13cm,width=9cm,angle=270]{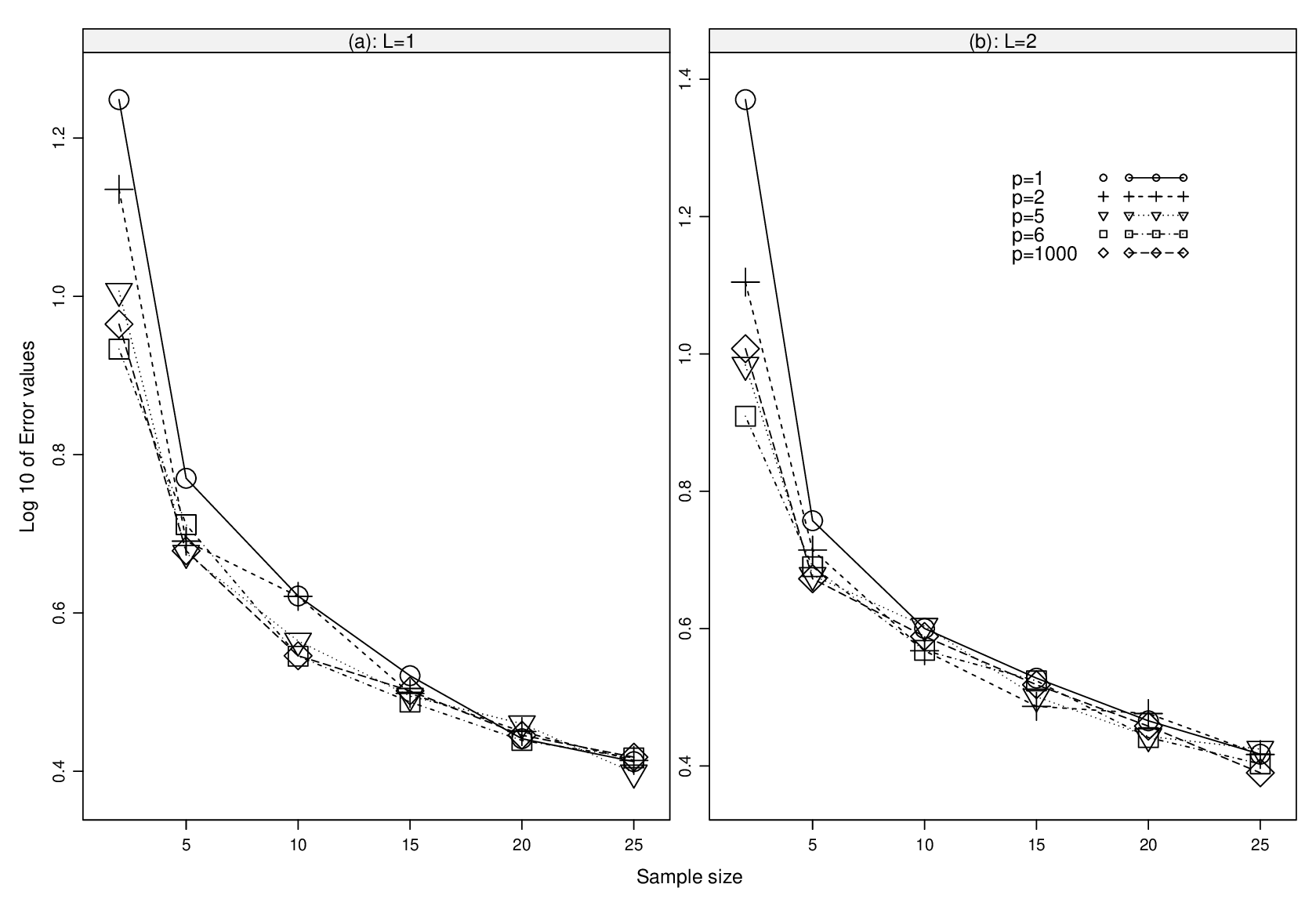}    
\end{center} 
\caption{Average of error values against the sample sizes ($N$) for $p$-balls (i.e., $\tilde{\bo{V}}$) and when $\bo{d=100}$. Panel (a) is obtained when $L=1$, while Panel (b) is associated with $L=2$.}     
 \label{fig:gradball100}         
\end{figure}

\begin{figure}[!hbp]           
\begin{center}
\includegraphics[height=13cm,width=9cm,angle=270]{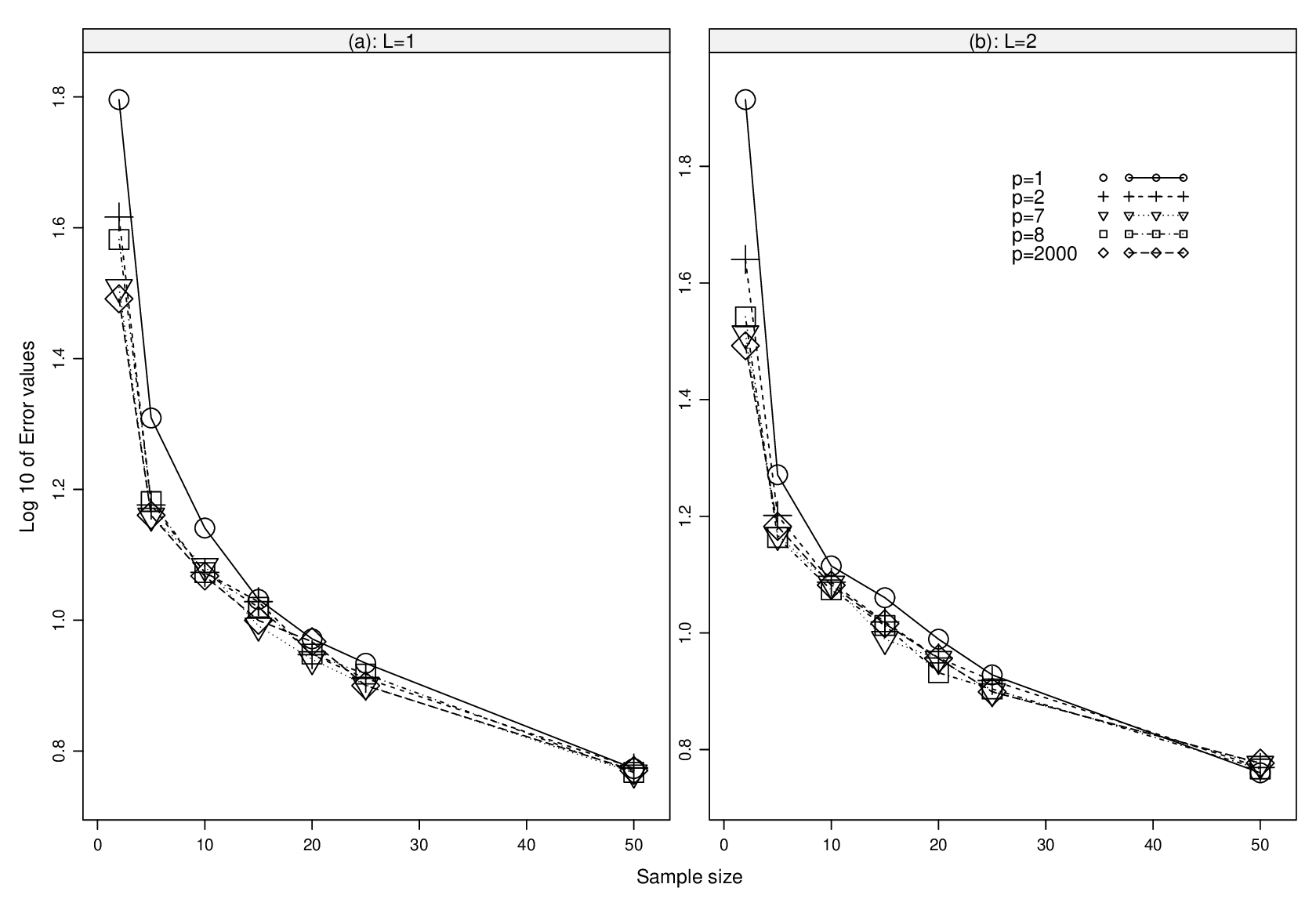}    
\end{center} 
\caption{Average of error values against the sample sizes ($N$) for $p$-spherical distributions (i.e., $\bo{V}$) when $\bo{d=1000}$. Panel (a) is obtained when $L=1$, while Panel (b) is associated with $L=2$.} 
 \label{fig:gradsp1000}         
\end{figure}

\begin{figure}[!hbp]         
\begin{center}
\includegraphics[height=13cm,width=9cm,angle=270]{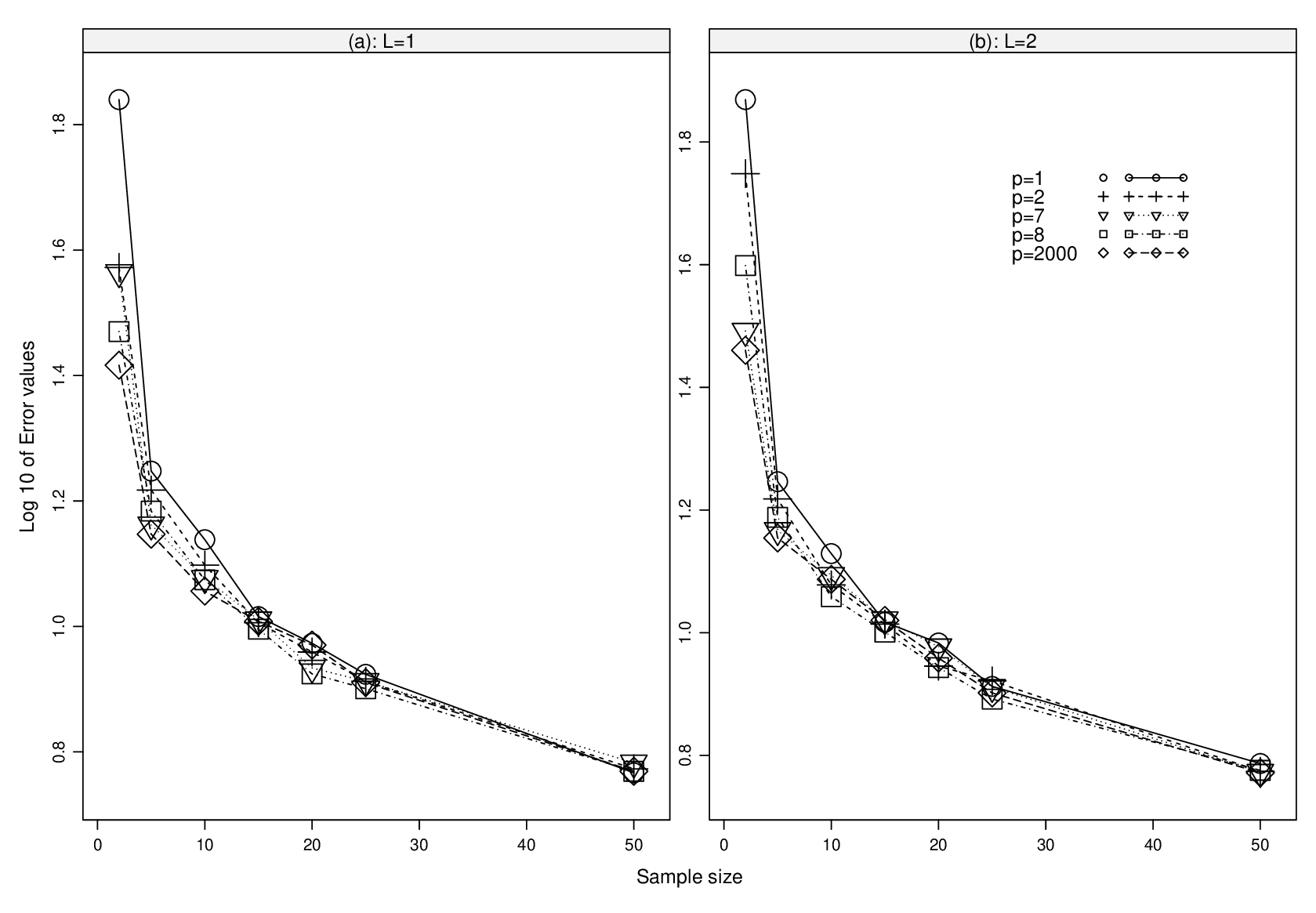}    
\end{center} 
\caption{Average of error values against the sample sizes ($N$) for $p$-balls (i.e., $\tilde{\bo{V}}$) and when $\bo{d=1000}$. Panel (a) is obtained when $L=1$, while Panel (b) is associated with $L=2$.}     
 \label{fig:gradball1000}            
\end{figure}

Based on Tables \ref{tab:fdm1}-\ref{tab:fdm4}, the proposed approach provides reasonable, accurate numerical results compared to other methods once the Gram-schmidt procedure is applied. Note that the error values are high without the Gram-schmidt procedure.  To assess the impact of $p$ on the error values
Figures \ref{fig:gradsp100}-\ref{fig:gradball1000} compare the errors for different values of $p$  when using the $p$-spherical distribution or $\tilde{\bo{V}}$. All of the error values are similar when  the Gram-schmidt procedure is applied, that is, when $N \geq d$ (not depicted here). Differences occur when $N \ll d$, and such differences are in favor of the derived $p$ (see Corollaries \ref{coro:optderi1}-\ref{coro:optderi2}). Thus, being able to generate empirical, uncorrelated random vectors equivalent to uniform distributions over the $p$-balls or spheres is necessary so as to improve such estimations. 
      
\subsection{A synthetic function} 
Consider the  synthetic function proposed in \cite{berahas22} and given by 
$$
\M_s(\bo{x}) := \sum_{k=1}^{d/2=100} \left[M_2 \sin(x_{2k-1}) + \cos(x_{2k}) \right] + \frac{M_1-M_2}{400} \bo{x}^\T \mathsf{1}_{200\times 200} \ \bo{x}  \, , 
$$     
with $M_1, M_2 \in \R$ and $\mathsf{1}_{200\times 200} \in \R^{200\times 200}$ the unit matrix (i.e., its entries are one). 
Its traditional gradient at $\bo{0}$ is $\nabla\M_s(\bo{0}) = M_2\left[1, 0,\, 1,\, 0\ldots, 1, 0 \right]^\T \in \R^{200}$. Tables \ref{tab:fdm5}-\ref{tab:fdm6} report the corresponding  error values. 
\begin{table}[htbp]  
\begin{center}   
\begin{tabular}{lcccc}   
\hline      
\hline          
 \textbf{ d=200} $G=\mathcal{I}$   &     \multicolumn{3}{c}{Number of total model evaluations (i.e., $LN$)} \\   
				&   $201$   &  $400$ & $400$  \\      			
 Methods &       &   &  \\      
\hline	
 FDM  (\cite{gilbert19})  &   -   &  - &   0.00005   \\   
\hline 
				 &   $L=1$   &   $L=1$  &  $L=2$   \\   
	  MC (\cite{lamboni24axioms}), $h=10^{-4}$  & 0.0049    &  0.0025   & 0.005  \\ 
	  This paper ($p=6$)  &   0.0049   &  0.0025  &     0.0049  \\           
\hline                                     
\hline                                
    \end{tabular}                
  \end{center}    
  \caption{Values of $Err$ for different approximations of $\nabla\M_s(\bo{0})$ using the  synthetic function associated with $\bo{M_1=2,\, M_2=1}$.}  
			 	\label{tab:fdm5}      
\end{table} 
\begin{table}[htbp]  
\begin{center}   
\begin{tabular}{lcccc}   
\hline      
\hline          
 \textbf{ d=200} $G=\mathcal{I}$   &     \multicolumn{3}{c}{Number of total model evaluations (i.e., $LN$)} \\   
				&   $201$   &  $400$ & $400$  \\      			
 Methods &       &   &  \\      
\hline	
 FDM  (\cite{gilbert19})  &   -   &  - &  0.049   \\   
\hline 
				 &   $L=1$   &   $L=1$  &  $L=2$   \\   
	  MC (\cite{lamboni24axioms}), $h=10^{-4}$  & 0.0065 & 0.0034   & 0.0049  \\ 
	  This paper ($p=6$)  &   0.0093   & 0.0027  &  0.0095\\           
\hline                                     
\hline                                
    \end{tabular}                 
  \end{center}    
  \caption{Values of $Err$ for different approximations of $\nabla\M_s(\bo{0})$ using the  synthetic function associated with $\bo{M_1=200,\, M_2=1/1000}$.}  
			 	\label{tab:fdm6}      
\end{table} 
It turns out that the proposed approach improves the results obtained in \cite{berahas22} using MC approaches. The FDM  performs well for this example when $M_1=2,\, M_2=1$ and vice versa when $M_1=200,\, M_2=1/1000$.       
                                                                      
\section{Conclusion} \label{sec:con}
In this paper, enhanced stochastic approximations of the gradients of functions evaluated at non-independent variables have been investigated. The proposed approach relies on a set of $L$ constraints and a class of $\ell_p$-spherical distributions or uniform distributions over the $p$-balls. The convergence analysis shows that the biases and the MSEs of the proposed estimators of dependent gradients (including the traditional gradient) of $2$-smooth functions do not suffer from the curse of dimensionality by properly choosing the values of $p$. Such results  outperform the best known results and allow for breaking down the course of dimensionality. \\  
 
Numerical results confirmed the reasonable accuracy of the proposed approach based on high-dimensional test cases.  To reach such a performance, the Gram-Schmidt algorithm is applied so as to obtain perfect, empirical and uncorrelated random vectors. Such a procedure limits the potential deployment of our algorithm to compute the gradients using lower model runs. It is then interesting to investigate  ways of generating perfect, empirical and uncorrelated random vectors equivalent to uniform distributions over the $p$-balls or on the $p$-spheres for a given sample size $N \ll d$. Moreover, noisy, non-smooth and stochastic functions are subjects for future investigations so as to extend the proposed approach. Also, one-point residual feeeback methods need extensions using the proposed estimators.     
              
\section*{Acknowledgments}     
  We would like to thank the three reviewers for their comments and remarks that  have helped improving this paper.           
\begin{appendices}       

\section{Proof of Corollary \ref{coro:parderord}} \label{app:coro:parderord}
It is shown in \cite{lamboni24axioms} (see also supplementary documents) that the bias  $B$ of the gradient satisfy      
\begin{eqnarray} 
 B & =&  \norme{ G^{-1}(\bo{x})  \esp \left[\nabla\M(\bo{x})-  \sum_{\ell=1}^{L} C_{\ell} \M\left(\bo{x} + \beta_\ell \boldsymbol{\hh}\bo{V} \right) \frac{\bo{V} \bo{h}^{-1}}{\sigma^2} \right]} \nonumber \\
 &\leq &   \sum_{\ell=1}^{L} \left| C_{\ell}\right|  \beta_\ell^2 M_2 \norme{ \left|G^{-1}(\bo{x})\right| \esp  \left[ \frac{\bo{V}^2}{\sigma^2} \norml{\boldsymbol{\hh}\bo{V}}\right]} \nonumber \, .
\end{eqnarray}        
Since $\hh_j=\hh$, $\bo{V} =R \bo{U}$, and using the properties of $\bo{U}$ provided in \cite{valle12} (Lemma 4.6), that is,   
$$
 \esp\left[ |U_i|^q \right] = \esp\left[ |U_1|^q \right] =\frac{\Gamma((q+1)/p) \Gamma(d/p)}{\Gamma(1/p) \Gamma(d/p+q/p) };
\qquad 
 \esp\left[ U_1^2 |U_2| \right] =  \frac{\Gamma(3/p) \Gamma(2/p) \Gamma(d/p)}{\Gamma^2(1/p) \Gamma(d/p+3/p) }
$$   
we can write    
\begin{eqnarray} 
 B  & \leq & \frac{ M_2 \hh \, \esp\left[ R^3 \right]}{\sigma^2} \norme{ \left|G^{-1}(\bo{x})\right| \esp  \left[ \bo{U}^2 \norml{\bo{U}}\right]}  \sum_{\ell=1}^{L} \left| C_{\ell}\right|  \beta_\ell^2   \nonumber \\
 & = & \frac{ M_2 \hh \, \esp\left[ R^3 \right]}{\sigma^2} \esp  \left[ |U_1|^3 + (d-1) U_1^2 |U_2|\right]
\norme{ \left|G^{-1}(\bo{x})\right| \indic_{\bu}}  \sum_{\ell=1}^{L} \left| C_{\ell}\right|  \beta_\ell^2 \nonumber \\
&=&   \frac{ M_2 \hh \, \esp\left[ R^3 \right]}{\sigma^2} \norme{ \left|G^{-1}(\bo{x})\right| \indic_{\bu}} 
\frac{\Gamma(d/p) \,\left[\Gamma(4/p) \Gamma(1/p) + (d-1)\Gamma(3/p) \Gamma(2/p)\right]}{\Gamma^2(1/p) \Gamma(d/p+3/p) } 
 \nonumber \, , 
\end{eqnarray}     
and the result holds.

\section{Proof of Lemma \ref{lam:qthm}} \label{app:lam:qthm}
One can see that $R_0 \sim \mathcal{U}(0,\, \xi)$ leads to $\esp\left[R_0^q\right] =  \frac{3^{q/2} \sigma^q}{q+1}
\left( \frac{\Gamma(1/p) \Gamma(d/p+2/p)}{\Gamma(3/p) \Gamma(d/p)}  \right)^{q/2}$.\\  
Point (i) relies on the approximation (\cite{gradshteyn07}, P. 904)     
$$
\frac{\Gamma(x+y)}{\Gamma(x)} = x^y \frac{\prod_{k=1}^{+\infty} \left(1+\frac{1}{x+k} \right)^y \exp\left(-y/(x+k) \right)}{\prod_{k=1}^{+\infty} \left(1+\frac{y}{x+k} \right) \exp\left(-y/(x+k) \right)}
\approx x^y \, ,      
$$ 
for higher values of $x$, corresponding to the choice $1\leq p\ll d$.\\
Point (ii) is the consequence of $\Gamma(1/p) \approx p$  for higher values of $p$.

\section{Proof of Corollary \ref{coro:dimfreup}} \label{app:coro:dimfreup}
 It is straightforward using Corollary \ref{coro:parderord} and Equation (\ref{eq:rmo}). Indeed, under $1 \leq p \ll d$,  and using $A_0 := K_{1,d,p} \frac{\esp\left[R_0^3\right]}{\sigma^2}$, we can write 
\begin{eqnarray}
A_{0} & = &  \frac{3\sqrt{3} \sigma^3}{4 \sigma^2} 
\frac{\Gamma(d/p) \,\left[\Gamma(4/p) \Gamma(1/p) + (d-1)\Gamma(3/p) \Gamma(2/p)\right]}{\Gamma^2(1/p) \Gamma(d/p+3/p) } 
\left( \frac{\Gamma(1/p) \Gamma(d/p+2/p)}{\Gamma(3/p) \Gamma(d/p)}  \right)^{3/2}  \nonumber \\
& := &  \frac{3\sqrt{3} \sigma}{4}
\frac{\left[\Gamma(4/p) \Gamma(1/p) + (d-1)\Gamma(3/p) \Gamma(2/p)\right]}{ \Gamma^{1/2}(d/p) \Gamma^{1/2}(1/p) \Gamma(d/p+3/p) }  
\left( \frac{\Gamma(d/p+2/p)}{\Gamma(3/p)}  \right)^{3/2}  \nonumber \\  
& \approx & 
\frac{\Gamma(d/p) \,\left[\Gamma(4/p) \Gamma(1/p) + (d-1)\Gamma(3/p) \Gamma(2/p)\right]}{\Gamma^2(1/p) \Gamma(d/p+3/p) } 
 \frac{3\sqrt{3}\sigma \, d^{3/p}}{(4) \, p^{3/p}}
\left( \frac{\Gamma(1/p)}{\Gamma(3/p)} \right)^{3/2}   \nonumber  \\
& = &  \frac{3\sqrt{3} \sigma}{4} 
\frac{\left[\Gamma(4/p) \Gamma(1/p) + (d-1)\Gamma(3/p) \Gamma(2/p)\right]}{\Gamma^{1/2}(1/p) \Gamma^{3/2}(3/p)}  
\nonumber  \, ,          
\end{eqnarray}
which yields the results.   \\ 
For  $1 \leq d \ll p$, we can see that $K_{1,d,p} \approx \frac{(d+3)(2d+1)}{12d}$, leading to 
$$
K_{1,d,p} \frac{\esp\left[R_0^3\right]}{\sigma^2} \approx \frac{(d+3)(2d+1)}{12d} \frac{3^3\sigma^3}{4\sigma^2} \left(\frac{d}{d+2}\right)^{3/2} =  \sigma \frac{9(d+3)(2d+1) d^{1/2}}{16(d+2)^{3/2}}   \, ,          
$$  
thanks to Equation (\ref{eq:rmo1}).

\section{Proof of Lemma \ref{lem:mosp}} \label{app:lem:mosp} 
Using the representation $\bo{V}= R \bo{U}$, the map $\bo{U} \mapsto g\left(R \bo{U} \right)$ is then a $L_0 R$-Lipschitz function w.r.t. $\bo{U}$. Denote with $\esp_U$ the expectation taking w.r.t. $U$, and $N_p$ the $p$-generalized standard normal variable (\cite{valle12}).  As $\bo{U}$ is a $p$-exponential concentrated random vector as well as $g(R\bo{U})$ (\cite{ledoux01,louart20}), taking the exponential inequality associated with $g(R\bo{U})$ yields 
\begin{eqnarray}    
A_1 & :=& \esp_{U}\left[\left|g\left(\bo{V} \right)- g\left(\bo{0} \right) \right|^q\right] = \int_{0}^{+\infty} \mathbb{P}\left[\left|g\left(\bo{V} \right)- g\left(\bo{0} \right)\right|^{q} \geq t \, | R\right] \, dt      \nonumber \\  
& = & \int_{0}^{+\infty} \mathbb{P}\left[\left|g\left(\bo{V} \right)- g\left(\bo{0} \right)\right| \geq t^{1/q} \, | R\right] \, dt      \nonumber \\  
& \stackrel{\cite{ledoux01,louart20}}{\leq} & C  \int_{0}^{+\infty} \exp\left(-\frac{c_0 d t^{p/q}}{L_0^pR^p} \right) \, dt = C q  \left(\frac{L_0^p R^p}{c_0 d p}\right)^{q/p} 
\int_{0}^{+\infty} y^{q-1}\exp\left(-\frac{y^p}{p} \right) \, dy   \nonumber \\     
& = &  \frac{C q}{2}  \left(\frac{L_0^p R^p}{c_0 d p}\right)^{q/p}  
\int_{-+\infty}^{+\infty} |y|^{q-1}\exp\left(-\frac{|y|^p}{p} \right) \, dy 
 =  \frac{C q}{2 \alpha_p}  \left(\frac{L_0^p R^p}{c_0 d p}\right)^{q/p}  \esp\left[ |N_p|^{q-1}\right] 
 \nonumber  \\        
& = &  \frac{C q}{2 \alpha_p}  \left(\frac{L_0^p R^p}{c_0 d p}\right)^{q/p}
\frac{p^{(q-1)/p} \Gamma(q/p)}{\Gamma(1/p)} =: K(q, p) \frac{L_0^q R^q}{d^{^{q/p}}}
  \nonumber  \, ,                
\end{eqnarray}       
because the $q$th-order moment of the $p$-generalized normal variable $N_p$ is given by $\esp\left[ |N_p|^{q}\right] = \frac{p^{q/p} \Gamma(q/p +1/p)}{\Gamma(1/p)}$ (see \cite{valle12}, Lemma 4.6). The quantity $\alpha_p := 0.5 p^{1-1/p}/ \Gamma(1/p)$ is the normalization constant. Thus,
 $$ 
K(q, p) :=  \frac{C q \Gamma(1/p)}{(c_0 p)^{q/p} p^{1-1/p}} \frac{p^{(q-1)/p} \Gamma(q/p)}{\Gamma(1/p)}
= \frac{q C  \Gamma(q/p)}{p(c_0)^{q/p}}    \, . 
$$         
The result holds by taking the expectation $\esp_{R}\left[A_1 \right]$ because $R$ and $\bo{U}$ are independent.       

\section{Proof of Corollary \ref{coro:mosp}} \label{app:coro:mosp}
 
Using Equation (\ref{eq:rmo}) and Lemma \ref{lem:mosp}, we can write 
\begin{eqnarray}
\esp\left[\left|g\left(\bo{V} \right)- g\left(\bo{0} \right)\right|^q\right]  &\lessapprox & 
\frac{q C  \Gamma(q/p)}{p c_0^{q/p}} \frac{L_0^q}{d^{^{q/p}}} \frac{ (3 \sigma^2)^{q/2}\, d^{q/p} }{(q+1) \, p^{q/p}} \left( \frac{\Gamma(1/p)}{\Gamma(3/p)} \right)^{q/2} \nonumber \\
&=& \frac{q C  \Gamma(q/p) (3 \sigma^2)^{q/2} L_0^q}{(q+1) p^{1+q/p} c_0^{q/p} } \left( \frac{\Gamma(1/p)}{\Gamma(3/p)} \right)^{q/2} \, , \nonumber      
\end{eqnarray}     
 and   
\begin{eqnarray}
\esp\left[R^q \left|g\left(\bo{V} \right)- g\left(\bo{0} \right)\right|^q\right]  &\lessapprox & 
\frac{q C  \Gamma(q/p) (3 \sigma^2)^{q} L_0^q \, d^{^{q/p}}}{(2q+1) p^{1+2q/p} c_0^{q/p} } \left( \frac{\Gamma(1/p)}{\Gamma(3/p)} \right)^{q} \, . \nonumber           
\end{eqnarray}      
Point (ii) is a consequence of Equation (\ref{eq:rmo1}) knowing that $\Gamma(1/p) \approx p$ when $p\gg 1$.

\section{Proof of Lemma \ref{lem:molpuni}} \label{app:lem:molpuni}
Firstly, as $
 \esp\left[ U_1^4\right] = \frac{\Gamma(5/p) \Gamma(d/p)}{\Gamma(d/p+4/p)} 
$ 
(see \cite{valle12}, Lemma 4.6 ), one can write 
$$    
 \esp\left[ \left|\left|  \bo{U} \right|\right|_2^4 \right]  \leq d \sum_{i=1}^d  \esp\left[U_i^4 \right] 
= d^2  \esp\left[U_1^4 \right] \stackrel{\cite{valle12}}{=} \frac{d^2 \Gamma(5/p) \Gamma(d/p)}{\Gamma(1/p) \Gamma(d/p+4/p)} \lessapprox  \frac{p^{4/p} \Gamma(5/p)}{\Gamma(1/p)}\, d^{2-4/p}\, .      
$$   
   
Secondly, using the concentration inequality for $\left|\left|  \bo{U} \right|\right|_2$ (see \cite{louart20}, Proposition 2.11), we have 
\begin{eqnarray} 
\esp\left[ \left|\left|  \bo{U} \right|\right|_2^q \right] &=&  \int_0^\infty \mathbb{P}\left[  \left|\left|  \bo{U} \right|\right|_2 > t^{1/q} \right] \, dt  \leq  C \int_0^\infty   \exp\left(- c_0' d\, t^{p/q} /d \right) \, dt   \nonumber \\ 
& \leq &   C \int_0^\infty   \exp\left(- c_0'\, t^{p/q} \right) \, dt = \frac{q C  \Gamma(q/p)}{p(c_0')^{q/p}} \, . \nonumber    
\end{eqnarray}          
  
\section{Proof of Theorem \ref{theo:mse}} \label{app:theo:mse}
Firstly, the  $MSE := \esp \left[\norme{ \widehat{grad(\M)}(\bo{x}) - grad(\M)(\bo{x})}^2 \right] $ is given by
$$  
MSE = \esp \left[\norme{\widehat{grad(\M)}(\bo{x}) - \esp\left[\widehat{grad(\M)}(\bo{x}) \right]}^2 + \norme{\esp\left[\widehat{grad(\M)}(\bo{x}) \right] - grad(\M)(\bo{x})}^2 \right] \, .    
$$      
Since the bias $\esp \left[ \norme{\esp\left[\widehat{grad(\M)}(\bo{x}) \right] - grad(\M)(\bo{x})}^2\right]$ has been already  derived, we focus ourselves on the second-order moment. By defining $g(\bo{V}) := \sum_{\ell=1}^{L} C_\ell \M\left(\bo{x}+ \beta_\ell  \hh \bo{V} \right)$ and knowing that $\sum_{\ell=1}^{L} C_\ell=0$, we can write      
$
 g(\bo{0}) = \sum_{\ell=1}^{L} C_\ell \M(\bo{x}) = 0    
$,   
leading to 
\begin{equation} \label{eq:equiv} 
Q(\bo{x}) :=
G^{-1}(\bo{x})  \frac{\bo{V}}{\sigma^2 \hh} 
\sum_{\ell=1}^{L} C_\ell \M(\bo{x}+ \beta_\ell  \hh\bo{V})
= G^{-1}(\bo{x})  \frac{\bo{V}}{\sigma^2 \hh}   
 \left[ g(\bo{V}) - g(\bo{0}) \right] \, .       
\end{equation}  
    
Using (\ref{eq:estgrad}), we can see that $ \esp\left[Q(\bo{x})\right] =\esp\left[\widehat{grad(\M)}(\bo{x}) \right]$. Bearing in mind the definition of the Euclidean norm and the variance, the centered second-order moment, that is, $\var_{grad} := \esp \left[\norme{\widehat{grad(\M)}(\bo{x}) - \esp\left[\widehat{grad(\M)}(\bo{x}) \right]}^2 \right]$ is given by   
\begin{eqnarray}       
\var_{grad} & \leq & \frac{1}{N}  \esp \left[\norme{G^{-1}(\bo{x})  \frac{\bo{V}}{\sigma^2 \hh} 
\sum_{\ell=1}^{L} C_\ell \M(\bo{x}+ \beta_\ell  \hh\bo{V}) - \esp\left[\widehat{grad(\M)}(\bo{x}) \right]}^2 \right] \nonumber \\   
& \leq & \frac{1}{N}  \esp \left[\norme{G^{-1}(\bo{x})  \frac{\bo{V}}{\sigma^2 \hh} 
\sum_{\ell=1}^{L} C_\ell \M(\bo{x}+ \beta_\ell \hh \bo{V})}^2 \right] \nonumber \\ 
& \stackrel{(\ref{eq:equiv})}{=} & \frac{1}{N}  \esp \left[\norme{G^{-1}(\bo{x})  \frac{\bo{V}}{\sigma^2 \hh}   
\left[ g(\bo{V}) - g(\bo{0}) \right]}^2 \right] \nonumber \\
& \leq & \frac{1}{N}  \esp \left[\norme{G^{-1}(\bo{x})  \frac{\bo{V}}{\sigma^2 \hh}}^2   
\left[ g(\bo{V}) - g(\bo{0}) \right]^2 \right] \nonumber \\ 
& \leq & \frac{1}{N}  \sqrt{\esp \left[\norme{G^{-1}(\bo{x})  \frac{\bo{V}}{\sigma^2 \hh}}^4 \right]}  \sqrt{
\esp \left[    
 \left\{    
g(\bo{V}) - g(\bo{0})
 \right\}^4
\right]}  \nonumber \\             
& \leq & \frac{1}{N \hh^2 \sigma^4}  \sqrt{ \esp\left[ R^4\right]} \sqrt{ \esp \left[\norme{G^{-1}(\bo{x})  \bo{U}}^4 \right]}  \sqrt{  
\esp \left[      
 \left\{    
g(\bo{V}) - g(\bo{0})      
 \right\}^4 
\right]}     \nonumber      \, .  
\end{eqnarray}    

Secondly, as $\M \in \mathcal{H}_{1}$, we have   
$    
\left| \M(\bo{x}+ \beta_\ell  \hh\bo{V}) -  \M(\bo{x}) \right| \leq M_{1} \norme{\beta_\ell  \hh\bo{V}}        
$, meaning that $\bo{V} \mapsto \M(\bo{x}+ \beta_\ell  \hh \bo{V})$ is a Lipschitz function with the constant $M_{1} |\beta_\ell| \hh$.  We can check that $\bo{V} \mapsto g(\bo{V}) = \sum_{\ell=1}^{L} C_\ell \M\left(\bo{x}+ \beta_\ell  \hh \bo{V} \right)$ is a Lipschitz function with the constant $L_0 := M_{1} \hh \sum_{\ell=1}^{L} | C_\ell \beta_\ell| = M_{1} \hh$. \\    
Finally, combining all these elements (i.e., Equation (\ref{eq:rmo}), Corollary \ref{coro:mosp} and Lemma \ref{lem:molpuni}) for $p \ll d$ yields 
\begin{eqnarray}       
\var_{grad} & \leq & \frac{1}{N \hh^2 \sigma^4}  \sqrt{ \frac{4 C  \Gamma(4/p)}{p c_0^{4/p}} \frac{M_1^4\hh^4}{d^{^{4/p}}}}    
 \frac{ (3 \sigma^2)^{2}\, d^{4/p} }{(5) \, p^{4/p}} \left( \frac{\Gamma(1/p)}{\Gamma(3/p)} \right)^{2}  \sqrt{ \esp \left[\norme{G^{-1}(\bo{x})  \bo{U}}^4 \right]}     \nonumber   \\ 
 & \leq & \frac{18 M_1^2}{5N}  \frac{C^{1/2}  \Gamma^{1/2}(4/p) \Gamma^2(1/p)}{c_0^{2/p}  p^{1/2+4/p} \Gamma^2(3/p)} d^{2/p}
\left|\left| G^{-1}(\bo{x}) \right|\right|_s^2  \sqrt{\esp \left[\norme{\bo{U}}^4 \right]}     \nonumber   \\ 
&= &  \frac{18 M_1^2}{5N}  \frac{C^{1/2}  \Gamma^{1/2}(4/p) \Gamma^{1/2}(5/p) \Gamma^{3/2}(1/p)}{c_0^{2/p}  p^{1/2+2/p} \Gamma^2(3/p)}   
\left|\left| G^{-1}(\bo{x}) \right|\right|_s^2  d    \nonumber  \, ,   
\end{eqnarray}   
using the second result of Lemma \ref{lem:molpuni}.  For the first result, we have    
$$
\var_{grad}  \leq \frac{36 M_1^2}{5N}  \frac{C_1^{1/2}  \Gamma(4/p) \Gamma^2(1/p)}{c_1^{2/p}  p^{1+4/p} \Gamma^2(3/p)} 
\left|\left| G^{-1}(\bo{x}) \right|\right|_s^2 d^{2/p} \, .     
$$   
Thirdly, under the assumption $1\leq d \ll p$, we can  write (thanks to Equation (\ref{eq:rmo1}), Corollary \ref{coro:mosp} and Lemma \ref{lem:molpuni}) 
\begin{eqnarray}       
\var_{grad} & \leq & \frac{1}{N \hh^2 \sigma^4} \frac{3^{4} \sigma^4 C' M_1^2 h^2}{5 c_0^{2/p}}\left( \frac{d}{d+2} \right)^{2}  d^{-2/p}    \nonumber   \\ 
 & \leq & \frac{81 C' M_1^2}{5 N} \left( \frac{d}{d+2} \right)^{2}  d^{-2/p}   \nonumber  \, ,   
\end{eqnarray}       
and the last result holds.

\end{appendices}                                   
         

\end{document}